\documentclass[12pt,a4paper,reqno]{amsart}
\usepackage{amsmath}
\usepackage{physics}
\usepackage{amssymb}
\numberwithin{equation}{section}
\usepackage{graphicx}
\usepackage{amssymb}
\usepackage{amsthm}
\usepackage{enumitem}
\usepackage{tikz}
\usetikzlibrary{shapes}
\usepackage{mathrsfs}
\usepackage{hyperref}
\usetikzlibrary{arrows}
\usepackage[utf8]{inputenc}
\usetikzlibrary{calc}
\numberwithin{equation}{section}
\usepackage{tikz-cd}
\usepackage[utf8]{inputenc}
\usepackage[T1]{fontenc}

\usepackage{tikz-cd}
\usepackage{txfonts}
\usepackage{lipsum} 
\let\oldpart\part
\renewcommand\part{\newpage\oldpart}

\usepackage[active]{srcltx}

\newtheorem{thm}{Theorem}[section]
\newtheorem*{theorem*}{Theorem}

\newtheorem{prop}[thm]{Proposition}
\newtheorem{lm}[thm]{Lemma}

\newtheorem{coro}[thm]{Corollary}

\newcommand{\name}[1]{{\bf #1}}

\newcommand{\fusion}[3]{{\binom{#3}{#1\;#2}}}

\newcommand\numberthis{\addtocounter{equation}{1}\tag{\theequation}}

\newtheorem{innercustomgeneric}{\customgenericname}
\providecommand{\customgenericname}{}
\newcommand{\newcustomtheorem}[2]{%
	\newenvironment{#1}[1]
	{%
		\renewcommand\customgenericname{#2}%
		\renewcommand\theinnercustomgeneric{##1}%
		\innercustomgeneric
	}
	{\endinnercustomgeneric}
}
\newcustomtheorem{customthm}{Theorem}

\theoremstyle{definition}
\newtheorem{example}[thm]{Example}
\newtheorem{df}[thm]{Definition}
\newtheorem{remark}[thm]{Remark}

\allowdisplaybreaks
\newcommand{\D}{\mathbb{D}}
\newcommand{\N}{\mathbb{N}}
\newcommand{\Z}{\mathbb{Z}}
\newcommand{\Q}{\mathbb{Q}}
\newcommand{\R}{\mathbb{R}}

\def \ra {\rightarrow}

\def \C {\mathbb{C}}

\def\la{\lambda}

\def \al{\alpha}
\def \si{\sigma}
\def \om{\omega}
\def \ga {\gamma}

\def \b {\beta}
\def \op {\oplus}

\def \ssq{\subseteq}
\def \vac {\mathbf{1}}

\def \g {\mathfrak{g}}
\def \h {\mathfrak{h}}

\def \End {\mathrm{End}}

\def\Id{\mathrm{Id}}
\def \wt {\mathrm{wt}}
\def\ep{\epsilon}

\def \bs {\backslash}

\def\ds{\dots}
\def\o{\otimes}

\def\spn{\mathrm{span}}

\def\fp{\mathfrak{p}}

\def\fn{\mathfrak{n}}
\def\fb{\mathfrak{b}}

\def\<{\langle}
\def\>{\rangle}
\def\Ind{\mathrm{Ind}}
\def\fb{\mathfrak{b}}
\def\Y{\mathcal{Y}}

\begin {document}

\title{Borel-type subalgebras of the lattice vertex operator algebra}
\author{Jianqi Liu}
\address{Department of Mathematics, University of Pennsylvania, Philadelphia, PA, 19104}
\email{jliu230@sas.upenn.edu}

\begin{abstract}
	In this paper, we introduce and study new classes of sub-vertex operator algebras of the lattice vertex operator algebras (VOAs), which we call the conic, Borel, and parabolic-type subVOAs. These CFT-type VOAs, which are not necessarily strongly finitely generated, satisfy properties similar to the usual Borel and parabolic subalgebras of a Lie algebra. For the lowest-rank nontrivial example of Borel-type subVOA $V_{B}$ of $V_{\Z\al}$, we explicitly determine its Zhu's algebra $A(V_B)$ in terms of generators and relations. 
\end{abstract}


\maketitle
\tableofcontents

\allowdisplaybreaks

\section{Introduction}

Conic, Borel, and parabolic-type subalgebras of the lattice vertex operator algebra (VOA) $V_L$ associated to a positive-definite even lattice $L$, defined here, are classes of CFT-type, irrational, $C_1$-cofinite, or non-$C_1$-cofinite VOAs that arise naturally in the Heisenberg module decomposition of a lattice VOA $V_{L}=\bigoplus_{\ga\in L} M_{\hat{\h}}(1,\ga)$ and give rise to solutions of the operator-form classical Yang-Baxter equations on the lattice VOAs. In this paper, we initiate a systematic investigation of these VOAs. In particular, we show that Borel and parabolic-type subVOAs share many properties similar to the usual Borel and parabolic subalgebras of a semisimple Lie algebra. To better understand the structure and representation theory of the Borel-type VOAs, we pick the smallest-rank nontrivial example of a Borel-type VOA $V_B=V_{\Z_{\geq 0}\al}$, a subVOA of the lattice VOA $V_{\Z\al}$, and explicitly determine its Zhu's algebra $A(V_B)$ in terms of generators and relations.

\subsection{Subalgebras of the lattice vertex operator algebra}

The lattice VOA $V_L$, introduced by Borcherds in \cite{B} and Frenkel, Lepowsky, and Meurman in \cite{FLM}, is the first example and plays a fundamental role in the theory of VOAs. The famous moonshine module vertex operator algebra $V^\natural$, which connects the $j$-invariant and the monster simple group $\mathbb{M}$, was constructed using the $\Z_2$-orbifold  and the simple current extension methods for the lattice VOA $V_{\Lambda}$ associated to the Leech lattice $\Lambda$ \cite{FLM}. Since then, the lattice VOA has been studied extensively: Dong and Lepowsky computed the fusion rules for $V_L$ in \cite{DL}. Dong classified the irreducible ($\theta$-twisted) modules over $V_L$ and proved the rationality of $V_L$ in \cite{D,D94}. 
The $\Z_2$-orbifold $V_L^+$ of the lattice VOA $V_L$, which was used to construct $V^\natural$, was studied in detail in a series of papers by Abe, Dong, Griess, Li, and Nagatomo \cite{AD04,ADL05,DG98,DN99(1),DN99}. Using the $\Z_2$ and $\Z_3$-orbifold constructions for the lattice VOAs, the Schellekens' conjecture on $71$ holomorphic VOAs with central charge $24$ \cite{S93} was solved recently by van Ekeren, Lam, Möller, Scheithauer, and Shimakura \cite{EMS,Lam11,LS15,LS16}.

Although the representation theory of lattice VOAs is well understood \cite{D,D94,DL}, their structural theory often leads to surprising results \cite{MS23}. On the other hand, it is evident that the structural theory of VOAs resembles semisimple Lie algebras. For instance, the Jacobi identity for VOAs is a formal variable generalization of the usual Jacobi identity \cite{FLM}; the vertex operator $Y$ satisfies the skew-symmetry axiom similar to a Lie bracket \cite{FHL}; and a CFT-type simple VOA such that $L(1)V_1=0$ admits a non-degenerate symmetric invariant bilinear form similar to a Cartan-Killing form \cite{FHL,L94}, etc. It is natural to expect the existence of substructures similar to Borel and parabolic-subalgebras of a semisimple Lie algebra in the lattice VOAs and others. Our definitions of the Borel-type subVOA $V_B$ and parabolic-type subVOA $V_P$, together with an auxiliary notion of conic-type subVOA $V_C$, of a lattice VOA $V_L$ is the first attempt in finding the analog of these substructures. 
We justify our terminologies by showing that they satisfy similar properties with the usual Borel and parabolic subalgebras of a semisimple Lie algebra in Section~\ref{Sec:3} and Sections~\ref{Sec:4}. These subVOAs of $V_L$ first appeared in solving the operator-form classical Yang-Baxter equation (CYBE) on VOAs by Bai, Guo, the author, and Wang in  \cite{BGL,BGLW}. We proved that Borel and parabolic-type subVOAs of a lattice VOA $V_L$ give rise to solutions of the CYBE on $V_L$. 

Besides the analogy with classical notions, our definition for the conic, Borel, and parabolic-type subVOAs also leads to large classes of CFT-type, $C_1$-cofinite or non-$C_1$-cofinite examples of VOAs. The $C_1$-cofiniteness, introduced by Li in \cite{L99}, is equivalent to the strongly finite-generation property for vertex algebras defined by Kac in \cite{K}. This finiteness condition for VOAs is satisfied by almost all the classical examples of CFT-type VOAs, and is essentially weaker than the more commonly acknowledged $C_2$-cofinite condition \cite{Z}. Although being a weak condition, $C_1$-cofiniteness is sufficient for many good properties of VOAs, like the convergence of iterated intertwining operators \cite{H05}, the Noetherianity of Zhu's algebra $A(V)$ \cite{Liu25}, the Morita equivalence between $A(V)$ and the mode transition algebra \cite{DGK23,DGK24}, and the Schur's Lemma for representations of VOAs over an arbitrary field \cite{YL23}, etc. It is expected that a $C_1$-cofinite VOA whose irreducible modules are also $C_1$-cofinite can give rise to a general conformal field theory that is not logarithmic, and the module category of such VOAs remains to be explored. Adding to the known examples of $C_1$-cofinite VOAs, we prove that the conic-type subVOAs of small rank lattice VOAs are $C_1$-cofinite, see Section~\ref{Sec:5.1}. On the other hand, due to the weakness of the $C_1$-cofinite condition, natural examples of CFT-type, non-$C_1$-cofinite VOAs are not easy to construct. As an application of our definitions, we find an example of a rank-two Borel-type subVOA $V_B$ that is not $C_1$-cofinite, see Section~\ref{Sec:5.2}. Examples arising from our definitions were also used to justify the necessity of the $C_1$-cofinite condition for the Noetherianity of Zhu's algebra $A(V)$ \cite{Liu25}. 

\subsection{The main theorems}
To state our main results, and describe how they are proved, we introduce some notation. Let $L$ be a positive-definite even lattice in an Euclidean space $E$. Our definitions are based on the key observation that the direct sum of the Heisenberg modules $V_M:=\bigoplus_{\ga\in M}M_{\hat{\h}}(1,\ga)$ associated to any (abelian) submonoid $M$ of the lattice $L$ is a CFT-type subVOA of $V_L$, whose vertex operator $Y_{V_M}$ is given by the intertwining operators among the Heisenberg modules $M_{\hat{\h}}(1,\ga)$ for $\ga\in M$. In this way, we can formulate plenty of new examples of subVOAs of $V_L$. 
We say that a submonoid $C\leq L$ is of conic-type if $C=\Z_{\geq 0}\al_1\op \ds\op \Z_{\geq 0}\al_r$, where $\{\al_1,\ds, \al_r\}$ is a basis of $L$; a submonoid $B\leq L$ is of Borel-type if $B\cup (-B)=L$, $B\cap (-B)=0$, and $B$ is contained in one side of a hyperplane in $E$; a submonoid $P\leq L$ is of parabolic-type if it contains a Borel-type submonoid, see Definition~\ref{def:submonoidsofaPDElattice}. The corresponding associated subVOAs $V_C$, $V_B$, and $V_P$ of the lattice VOA $V_L$ are called the conic, Borel, and parabolic-type subVOAs, respectively. 

Although there are many classes of the Borel-type subVOA of a given lattice VOA $V_L$, and their properties are very different, we can still justify our terminologies in terms of analogies and relations with the usual Borel and parabolic subalgebras of a semisimple Lie algebra in Section~\ref{Sec:4}. First, we prove that there is a correspondence between the degree-one sub-Lie algebras $(V_B)_1$ and $(V_P)_1$ of the (reductive) Lie algebra  $\g=(V_L)_1$ \cite{DM04} and the usual Borel and parabolic subalgebras of $\g$, see Propositions~\ref{prop:Borelfirstlevel1} and \ref{prop:Borelfirstlevel2}. Then we prove that conic-type subVOAs $V_C$ of the ADE-type root lattices $Q$ are conjugate under $\mathrm{Aut}(V_Q)$, see Proposition~\ref{prop2.10}. Next, we prove that conic, Borel, and parabolic-type subVOAs are not simple nor rational, see Proposition~\ref{prop:nonsimpleofconic} and Proposition~\ref{prop2.7}, which are parallel to the facts that the Borel and parabolic subalgebras of a Lie algebra are not simple nor semisimple. Then, we introduce a notion of ``normalizer'' $N_{W}(V)$ of a sub-vertex algebra $W\leq V$, which is similar to the notion of ``centralizer'' or commutant $\mathrm{Com}_V(W)$ introduced by Frenkel and Zhu in \cite{FZ} and is a generalization of the usual normalizer of a subalgebra in a Lie algebra. We prove that $N_{V_L}(V_P)=V_P$ for any parabolic-type subVOA $V_P\leq V_L$, see Proposition~\ref{prop2.13}, which is parallel to the facts that $N_{\g}(\fp)=\fp$ and $N_G(P)=P$, where $\fp\leq \g$ is a parabolic subalgebra of a semisimple Lie algebra $\g$, and $P\leq G$ is a parabolic subgroup of a linear algebraic group $G$ \cite{BT,Hum1}. 

The $C_1$-cofiniteness of these subVOAs, which has no analogs in the classical Lie theory, are discussed in Section~\ref{Sec:5}. Since these VOAs are irrational, it is well indicated that they are not $C_2$-cofinite. 
The following is our first main theorem, which concerns the $C_1$-cofiniteness of small-rank conic-type VOAs, see Lemma~\ref{lm2.14} and Theorem~\ref{thm:ranktwocoincC1}:
\begin{customthm}{A}\label{thm:A}
	The conic-type subVOA $V_C=V_{\Z_{\geq 0}\al}$ of the rank-one lattice VOA $V_{\Z\al}$ is $C_1$-cofinite.	
	
	Moreover, let  $L=\Z\al\op \Z\b$ be a rank-two positive-definite even lattice such that $(\al|\b)=-n$, $(\b|\b)=2k$ and $(\al|\al)=2\ell$, with $n\geq 1$, and $k\geq \ell\geq 1$. Suppose 
	$n^2+\ell^2-4\ell k\leq 0,$ then the conic-type subVOA $V_C=V_{\Z_{\geq 0}\al\op \Z_{\geq 0}\b}$ is $C_1$-cofinite. 
\end{customthm}

Theorem \ref{thm:A} leads to new examples of CFT-type $C_1$-cofinite VOAs.
A similar argument also proves the $C_1$-cofiniteness of a typical parabolic-type subVOA $V_P=V_{\Z\al\op \Z_{\geq 0}\b}$ of the lattice VOA $V_{A_2}$, see Proposition~\ref{prop:C1ofVP}. However, there is no uniform result on the $C_1$-cofiniteness of Borel-type subVOAs. Some of them are $C_1$-cofinite, while others are not. For instance, we prove that the rank-one Borel-type subVOA $V_{\Z_{\geq 0}\al}$ of $V_{\Z\al}$ is $C_1$-cofinite, while a rank-two Borel-type subVOA $V_B$ of the lattice VOA $V_{A_2}$, where $B=\{ m\al+n\b: m,n\in \Z, n>0\}+ \Z_{\geq 0}\al\leq A_2$, is not $C_1$-cofinite, see Proposition~\ref{prop:VBnonC1}.

To better understand the representation theory of our new examples of CFT-type VOAs, we focus on the smallest-rank nontrivial example of the Borel-type subVOA $V_{B}=V_{\Z_{\geq 0}\al}$, which is the only possible Borel-type subVOA of the rank-one lattice VOA $V_{\Z\al}$ by our definition. We determine the concrete structure of its Zhu algebra $A(V_B)$ in Section~\ref{sec3}. By construction, $A(V)=V/O(V)$, where $O(V)\subset V$ is spanned by elements of the form $a\circ b=\sum_{j\geq 0}\binom{\wt a}{j} a_{j-2}b$. Zhu proved that $A(V)$ is an associative algebra attached to a VOA $V$ such that irreducible $V$-modules are in one-to-one correspondence with irreducible $A(V)$-modules \cite{Z}. Using a sequence of inductions on the length of the spanning elements $u=\al(-n_1)\ds \al(-n_r)e^{m\al}$ of $V_B$, together with a detailed analysis for the relations between lattice vertex operators restricted to $V_B$ \cite{FLM}, which are carried out by Lemma~\ref{lm3.4}, Proposition~\ref{prop3.5}, Proposition~\ref{lm3.6}, and Proposition~\ref{prop3.7}, we provide a concrete description for $O(V_B)$ and show that $A(V_B)$ is a one-dimensional non-abelian nilpotent extension of the polynomial ring $\C[x]$. The following is our second main theorem, see Theorem~\ref{thm3.8} and Corollary~\ref{coro3.9}:

\begin{customthm}{B}\label{main:B}
	Let $V_B=V_{\Z_{\geq 0}\al}$ be the Borel-type subVOA of $V_{\Z\al}$, with $(\al|\al)=2N$. Then 
	\begin{equation}
		\begin{split}O(V_B)=\spn \bigg\{ \al(-n-2)u+\al(-n-1)u,\ \al(-1)v+v,\ M_{\hat{\h}}(1,k\al):&
			\\n\geq 0,\  u\in V_B, \ v\in \bigoplus_{m\geq 1}M_{\hat{\h}}(1,m\al),\ k\geq 2\bigg\}.
		\end{split}
	\end{equation} 
	In particular, $A(V_B)\cong \C[x]\op \C y$, with $y^2=0,yx=-Ny$ and $xy=Ny$. 
\end{customthm}

As applications of Theorem~\ref{main:B}, we prove that the irreducible modules of $V_B$ are in one-to-one correspondence with irreducible modules over the rank-one Heisenberg VOA $M_{\widehat{\C\al}}(1,0)$, and the fusion rules among irreducible modules over $V_B$ are the same as the fusion rules among irreducible modules over $M_{\widehat{\C\al}}(1,0)$, see Theorem~\ref{thm3.14}. This result is also parallel to the semisimple Lie algebra case, namely, the irreducible modules over a Borel subalgebra $\fb=\h\op \fn_+$ of $\g$ are the same as irreducible modules over the Cartan part $\h$ on which $\fn_+$ acts as $0$.

This paper is organized as follows: we first introduce the intermediate notions of conic, Borel, and parabolic-type submonoids of $L$ and discuss their basic properties in Section~\ref{Sec2}. Then we introduce the notions of conic, Borel, and parabolic-type subVOAs of the lattice VOA $V_L$ and give small-rank examples of these VOAs in Section~\ref{Sec:3}. In Section~\ref{Sec:4}, we prove the basic properties of these subVOAs which are analogous to the properties of the usual Borel and parabolic subalgebras. In Section~\ref{Sec:5}, we discuss the $C_1$-cofiniteness of these subVOAs and prove Theorem~\ref{thm:A}.  In Section~\ref{sec3}, we focus on the rank-one Borel-type subVOA $V_B=V_{\Z_{\geq 0}\al}$ of the rank-one lattice VOA $V_{\Z\al}$ and determine its Zhu's algebra $A(V_B)$ in terms of generators and relations.

\subsection{ Conventions}\label{sec:convention}
Throughout this paper, we adopt the following symbolic conventions: 
\begin{itemize}
	\item  All vector spaces are defined over $\C$, the complex number field. 
	\item $\N$ represents the set of all natural numbers, including $0$.
	\item  Unless we state otherwise, the letter $L$ represents a positive-definite even lattice in a Euclidean space $(E,(\cdot| \cdot))$.
	\item Let $E$ be a Euclidean space $E$. $O(E)=\{\si\in \mathrm{GL}(E):(\si(\al)|\si(\b))=(\al|\b),\ \forall \al,\b\in E\}$ is the group of isometries of $E$ \cite{CS88}.
	\item A VOA is said to be of {\em CFT-type} if $V=V_0\op V_+$. where $V_0=\C\vac$ and $V_+=\bigoplus_{n=1}^\infty V_n$  \cite{FHL,LL}.  
	\item Let $V$ be a VOA. $V'$ is the contragredient module of $V$ \cite{FHL}.
	\item $N_G(P)$ (resp. $N_\g(\fp)$) is the normalizer of the subgroup $P$ (resp. subalgebra $\fp$) in a group $G$ (resp. Lie algebra $\g$).
	\item $A_1$ and $A_2$ represent the root lattices of type $A_1$ and $A_2$, respectively \cite{Hum2}. 
\end{itemize}

\section{Submonoids of an even lattice}\label{Sec2}
Let $L$ be a lattice in $E$. By definition, there exists a $\R$-basis $\{\al_1,\ds, \al_r\}$ of $E$ such that $L=\Z\al_1\op\ds\op \Z\al_r$. Moreover, $(\al|\b)\in \Z$ for all $\al,\b\in L$, and $(\al|\al)\in 2\Z$ for all $\al\in L$.

\subsection{Definition of conic, Borel, and parabolic-type submonoids} We set some notation first. Let $\ga$ be a nonzero vector in the Euclidean space $E$. Let $P(\ga)$ be the hyperplane passing though the origin that is perpendicular to $\ga$, and let $P^+(\ga)$ (resp. $P^{\geq 0}(\ga)$) be the positive (resp. non-negative) side of $P(\ga)$. That is,  
\begin{align*}
	P(\ga):&=\{ v\in E: (\ga|v)=0\},\\
	P^+(\ga):&=\{ v\in E: (\ga|v)>0\},\numberthis\label{eq:hyperplanes} \\
	P^{\geq 0}(\ga):&=P(\ga)\sqcup P^+(\ga)= \{ v\in E: (\ga|v)\geq 0\}. 
\end{align*}

\subsubsection{Definition and first properties}
Observe $L$ is an abelian monoid with the product given by addition and the identity element given by $0$. An \name{(abelian) submonoid} of $L$ is a subset $M\subset L$ containing $0\in M$ that is closed under the addition of $L$.  An \name{(abelian) sub-semigroup} of $L$ is a subset $S\subset L$ that is closed under addition of $L$. We introduce the following notions: 

\begin{df}\label{def:submonoidsofaPDElattice}
	Let $L$ be a positive definite even lattice in a Euclidean space $(E,(\cdot| \cdot))$.
	\begin{enumerate}
		\item  A submonoid $C\leq L$ is said to be of {\bf conic-type} if there exists a basis $\{\al_1,\ds ,\al_r\}$ of $L$ such that $C=\Z_{\geq 0} \al_1\op \ds\op \Z_{\geq 0}\al_{r}$. 
		\item 	A submonoid $B\leq L$ is said to be of {\bf Borel-type} if it satisfies 
		\begin{enumerate}
			\item $B\cup (-B)=L$,
			\item $ B\cap (-B)=\{0\}$,
			\item There exists a hyperplane $P(\ga)\subset E$ such that $B\subset P^{\geq 0}(\ga)$. 
		\end{enumerate}
		\item A submonoid $P\leq L$ is said to be of \name{parabolic-type} if there exists a Borel-type submonoid $B\leq L$ such that $B\ssq P$. 
	\end{enumerate}
\end{df}
Recall the automorphism group of the lattice $L$ is defined by
\begin{equation}\label{eq:defofAutL}
	\mathrm{Aut}(L):=\{ \sigma\in \mathrm{GL}(L):(\si \al|\si\b)=(\al|\b),\ \forall \al,\b\in L\}=\mathrm{GL}(L)\cap O(E),
\end{equation}
which is a finite group \cite{CS88}. The following properties of these submonoids are immediate consequences of our definition. 
\begin{prop}\label{prop:propertiesofBorel}
	Conic, Borel, and parabolic-type submonoids of $L$ satisfy:
	\begin{enumerate}
		\item Every Borel-type submonoid is of parabolic-type. Every Borel-type submonoid contains a conic-type submonoid;
		\item The sets of conic, Borel, or parabolic-type submonoids of $L$ are invariant under the actions of $\mathrm{Aut}(L)$;
		\item The hyperplane $P(\ga)$ in the Definition~\ref{def:submonoidsofaPDElattice} of Borel-type submonoid is unique. 
	\end{enumerate}
\end{prop}
\begin{proof}
	The first claim  in $(1)$ is clear. Let $B$ be a Borel-type submonoid, and let $\{\al_1,\ds ,\al_r\}$ be a basis of $L$. Since $B\cup (-B)=L$, by replacing $\al_i$ with $-\al_i$, which does not change the fact that $\{ \al_1,\ds,\al_r\}$ is a basis of $L$, we may assume $\{\al_1,\ds ,\al_r\}\ssq B$, then $C=\Z_{\geq 0} \al_1\op \ds\op \Z_{\geq 0}\al_{r}$ is a conic-type submonoid that is contained in $ B$. 
	
	(2) Let $\si\in \mathrm{Aut}(L)$. Since $\si$ preserves addition, $\si(M)\leq L$ is a submonoid for any submonoid $M\leq L$. Moreover, since $\mathrm{Aut}(L)\subset \mathrm{SL}(L)$, $\{ \si(\al_1),\ds ,\si(\al_r)\}$ is another $\Z$-basis of $L$. Hence $\si (C)=\Z_{\geq 0} \si(\al_1)\op\ds\op \Z_{\geq 0}\si(\al_r)$ is another conic-type submonoid. Let $B\leq L$ be a Borel-type submonoid, we have $L=\si(L)=\si(B)\cup (-\si(B))$ and $\si(B)\cap (-\si(B))=\si(B\cap(-B))=\si(0)=0$. Furthermore, since $\si\in O(E)$ \eqref{eq:defofAutL}, it is easy to see that $\si (P^{\geq 0}(\ga))=P^{\geq 0}(\si(\ga))$, with $\si(\ga)\neq 0$. Thus, $\si(B)\subset P^{\geq 0}(\si(\ga))$, and $\si(B)$ is a Borel-type submonoid. Finally, it is clear that $\si(P)$ is a parabolic-type submonoid since $\si(B)\ssq \si(P)$. 
	
	(3) Let $B$ be a Borel-type submonoid such that $B\subset P^{\geq 0}(\ga)$. First, we claim that $L\cap P^{>0}(\ga)\ssq B$ and $L\cap P^{<0}(\ga)\ssq -B$. Indeed, for any $\al\in L\cap P^{>0}(\ga)$, since $\al\neq 0$ and the subset $B$ satisfies $B\cup (-B)=L$ and $B\cap (-B)=0$, we have $\al\in B$ or $\al\in -B$. If $\al\in -B$ then $(-\al|\ga)\geq 0$, which contradicts $(\al|\ga)>0$. Hence $\al\in B$ and so $L\cap P^{>0}(\ga)\ssq B$. 
	
	Suppose there exists another hyperplane $P(\ga')$, with $\ga'$ not being parallel to $\ga$, that defines the same Borel-type submonoid $B$. i.e., $B\cup (-B)=L$, $B\cap (-B)=\{0\}$, and $B\subset P^{\geq 0}(\ga')$. Assume $L=\Z\al_1\op\ds\op \Z\al_r$. Since $(\cdot |\cdot)$ is positive-definite on $E$, we may choose $\b_1,\ds,\b_r\in E$ such that $(\al_i|\b_j)=\delta_{i,j}$ for all $1\leq i,j\leq r$. Then $\{\b_1,\ds,\b_r\}$ is a $\R$-basis of $E$.
	
	Write $\ga= a_1\b_1+\ds +a_r\b_r$ and $\ga'=a_1'\b_1+\ds+a_r'\b_r$, where $\vec{a}=(a_1,\ds ,a_r)$ and $\vec{a'}=(a_1',\ds,a_r')$ are two non-parallel vectors in $\R^r\bs\{0\}$. Then 
	$$\D:=\{(\la_1,\ds,\la_r)\in \R^r: \la_1a_1+\ds+\la_ra_r<0,\ \la_1a_1'+\ds+\la_ra_r'>0 \}$$
	is a nonempty open subset in $\R^r$. Since $\Q^r\bs\{0\}$ is dense in $\R^r$, we have $\D\cap (\Q^r\bs\{0\})\neq \emptyset$. Therefore, by multiplying the common denominator, we can find a vector $(\la_1,\ds,\la_r)\in \Z^r\bs \{0\}$, which gives rise to a nonzero element $\al=\sum_{i=1}^r \la_i \al_i\in L$, such that 
	\begin{align*}
		(\al|\ga)&=\sum_{i,j}(\la_i\al_i|a_j\b_j)=\la_1a_1+\ds+\la_ra_r<0,\\
		(\al|\ga')&=\sum_{i,j}(\la_i\al_i|a'_j\b_j)=\la_1a_1'+\ds+\la_ra_r'>0.
	\end{align*}
	By the first inequality and our claim, we have $\al\in L\cap P^{<0}(\ga)\ssq -B$. While the second inequality indicates $\al\in L\cap P^{>0}(\ga')\ssq B$. Then $\la\in B\cap (-B)=0$, which is a contradiction.  Hence the hyperplane that defines Borel-type submonoid is unique. 
\end{proof}
\subsubsection{Existence of Borel-type submonoid}
Next, we show that Borel-type submonoids must exist. Hence the conic-type and parabolic-type submonoids also exist.  

\begin{lm}\label{lm:hyperplane}
	Let $r$ be a positive integer. There exists $\la_1,\ds,\la_r\in \R\bs\{0\}$, such that 
	\begin{equation}\label{eq:nonrationalcombo}
		m_1\la_1+m_2\la_2+\ds+m_r\la_r\notin \Q,  
	\end{equation} 
	for any nonzero rational vector $(m_1,\ds , m_r)\in \Q^r\bs\{0\}$. 
\end{lm}
\begin{proof}
	We use induction on $r$ to show the existence of $\la_1,\ds,\la_r$. When $r=1$, we may choose $\la_1\in \R\bs\Q$. Then clearly $m_1 \la_1\notin \Q$ for any nonzero rational number $m_1$. Assume there exists $\la_1,\ds,\la_{r-1}\in \R\bs\{0\}$ such that \eqref{eq:nonrationalcombo} holds for $r-1$. 
	
	Note that the $\Q$-vector space $\Q\la_1+\Q\la_2+\ds+\Q\la_{r-1}+\Q$ is a countable subset of $\R$ as the countable union of countable sets is countable. But $\R$ is uncountable. Hence there exists some element $\la_r\in \R\bs( \Q\la_1+\Q\la_2+\ds+\Q\la_{r-1}+\Q)$.  Now let $(m_1,\ds , m_r)\in \Q^r\bs\{0\}$. If $m_r=0$, then by the induction hypothesis, we have $m_1\la_1+m_2\la_2+\ds+m_{r-1}\la_{r-1}+0\la_r\notin \Q$ since not all $m_1,\ds,m_{r-1}$ are $0$. Assume $m_{r}\neq 0$. Since $\la_r\notin \Q\la_1+\Q\la_2+\ds+\Q\la_{r-1}+\Q$, we have 
	$$\frac{m_1}{m_r}\la_1+\frac{m_2}{m_r}\la_2+\ds +\frac{m_{r-1}}{m_r}\la_{r-1}+\la_r\notin \Q.$$
	Hence $m_1\la_1+m_2\la_2+\ds+m_r\la_r\notin \Q$. i.e., \eqref{eq:nonrationalcombo} holds for $r$. 
\end{proof}

\begin{prop}\label{prop:existenceofBorel}
	Every positive-definite even lattice $L$ contains a Borel-type submonoid. 
\end{prop}
\begin{proof}
	We claim that there exists a hyperplane $P(\ga)=\{ v\in E: (\ga|v)=0\}$ passing through the origin of $E$, such that $P(\ga)$ does not intersect any other lattice point of $L$. i.e., $P(\ga)\cap L=\{0\}$.
	
	Indeed, assume $L=\Z\al_1\op\ds\op \Z\al_r$. Similar to Proposition~\ref{prop:propertiesofBorel}, we choose a basis $\{\b_1,\ds,\b_r\} $ of $E$ such that $(\al_i|\b_j)=\delta_{i,j}$ for all $1\leq i,j\leq r$. Now choose $\la_1,\ds,\la_r\in \R\bs\{0\}$ such that \eqref{eq:nonrationalcombo} holds. Let $\ga:=\la_1\b_1+\ds+\la_r\b_r\in E$. Then for any $v=m_1\al_1+\ds+m_r\al_r\in L\bs\{0\}$, where $m_1,\ds, m_r\in \Z\subset \Q$ are not all $0$, we have 
	$$(\ga|v)=(\la_1\b_1+\ds+\la_r\b_r|m_1\al_1+\ds+m_r\al_r)=m_1\la_1+m_2\la_2+\ds+m_r\la_r\neq 0,$$
	in view of \eqref{eq:nonrationalcombo}. Hence $P(\ga)\cap L=\{0\}$. 
	
	Now let $P^+(\ga):=\{ v\in E: (\ga|v)>0\}$ be the ``positive'' half of the Euclidean space $E$ divided by $P(\ga)$, and let $B:=(L\cap P^+(\ga))\cup\{0\}$. Clearly, $B\leq L$ is a submonoid and $B\subset P^{\geq 0}(\ga)$. Moreover, since $E=P^+(\ga)\sqcup (-P^+(\ga))\sqcup P(\ga)$ and $L\cap P(\ga)=\{0\}$, we have $L=B\cup (-B)$ and $B\cap (-B)=\{0\}$. Thus, $B$ is a Borel-type submonoid.  
\end{proof}

In fact, the method in the proof of Proposition~\ref{prop:existenceofBorel} is {\bf not} the only way to construct Borel-type submonoids, as we will see in the following examples.
\begin{example}\label{ex:rankoneBorel}
	Let $L=\Z\al$ be a rank-one even lattice. Then $B=\Z_{\geq 0}\al$ is a submonoid of both conic and Borel-type. It is easy to see that $\Z_{\geq 0}\al$ and $\Z_{\leq 0} \al$ are the only possible Borel-type submonoids of $L$. In fact, these are the only possible parabolic-type submonoids as well. 
\end{example}

\begin{example}\label{ex:ranktwoparabolic}
	Let $L=A_2=\Z\al\op \Z\b$ be the type-$A_2$ root lattice \cite{FH91,Hum2}, where 
	$\{\al,\b\}$ is a set of simple roots with $(\al|\al)=(\b|\b)=2$ and $(\al|\b)=-1$. See   Figure~\ref{fig1}. Let 
	\begin{equation}\label{eq:exampleofB}
		B:=\{ m\al+n\b: m,n\in \Z, n>0\}+ \Z_{\geq 0}\al. 
	\end{equation}
	Clearly, $B\leq A_2$ is a submonoid. In Figure~\ref{fig1}, the red dots represent nonzero elements in $B$, while the blue dots represent nonzero elements in $-B$. It is clear from this graph that $A_2=B\cup (-B)$, $B\cap (-B)=\{0\}$, and $B$ is contained in $P^{\geq 0}(2\b+\al)$. Thus, $B$ is a Borel-type submonoid in $A_2$. Note that $B\bs\{0\}$ is not contained in the strictly positive side of any hyperplane passing through the origin, as opposed to the example constructed in Proposition~\ref{prop:existenceofBorel}. We call $B'=B\bs\{0\}$ the \name{Borel-type sub-semigroup} of $A_2$. 
	
	Let $P:=\Z\al\op \Z_{\geq 0}\b$. Then $P$ is a parabolic-type submonoid of $A_2$ since it contains the submonoid $B$ in \eqref{eq:exampleofB}. Note that $P$ is {\bf not} of Borel-type since $P\cap (-P)=\Z\al \neq \{0\}$. 
	
	\begin{figure}
		\centering
		\begin{tikzpicture}
			\coordinate (Origin)   at (0,0);
			\coordinate (XAxisMin) at (-5,0);
			\coordinate (XAxisMax) at (5,0);
			\coordinate (YAxisMin) at (0,0);
			\coordinate (YAxisMax) at (0,5);
			
			\clip (-5.2,-3.8) rectangle (5.2,3.7); 
			\begin{scope} 
				\pgftransformcm{1}{0}{1/2}{sqrt(3)/2}{\pgfpoint{0cm}{0cm}} 
				\coordinate (Bone) at (0,2);
				\coordinate (Btwo) at (2,-2);
				\draw[style=help lines,dashed] (-7,-6) grid[step=2cm] (6,6);
				\foreach \x in {-4,-3,...,4}{
					\foreach \y in {1,...,4}{
						\coordinate (Dot\x\y) at (2*\x,2*\y);
						\node[draw,red,circle,inner sep=2pt,fill] at (Dot\x\y) {};
					}
				}
				\foreach \x in {1,...,4}{
					\foreach \y in {0}{
						\coordinate (Dot\x\y) at (2*\x,2*\y);
						\node[draw,red,circle,inner sep=2pt,fill] at (Dot\x\y) {};
					}
				}
				\foreach \x in {-1,...,-4}{
					\foreach \y in {0}{
						\coordinate (Dot\x\y) at (2*\x,2*\y);
						\node[draw,blue,circle,inner sep=2pt,fill] at (Dot\x\y) {};
					}
				}
				
				\foreach \x in {-4,-3,...,4}{
					\foreach \y in {-1,...,-4}{
						\coordinate (Dot\x\y) at (2*\x,2*\y);
						\node[draw,blue,circle,inner sep=2pt,fill] at (Dot\x\y) {};
					}
				}
				\foreach \x in {0}{
					\foreach \y in {0}{
						\coordinate (Dot\x\y) at (2*\x,2*\y);
						\node[draw,circle,inner sep=2pt,fill] at (Dot\x\y) {};
						\node [xshift=0cm, yshift=-0.4cm]{$0$};
					}
				}
				
				\draw [thick,-latex,red] (Origin) 
				-- (Bone) node [above right]  {$\al+\b$};
				\draw [thick,-latex,blue] (Origin)
				-- (Btwo) node [below right]  {$-\b$};
				\draw [thick,-latex,red] (Origin) node [xshift=0cm, yshift=-0.3cm]{}
				-- ($(Bone)+(Btwo)$) node [below right] {$\al$};
				\draw [thick,-latex,blue] (Origin)
				-- ($-1*(Bone)-1*(Btwo)$) node [below left] {$-\al$};
				\draw [thick,-latex,red] (Origin)
				-- ($-1*(Btwo)$) coordinate (B3) node [above left] {$\b$};
				\draw [thick,-latex,blue] (Origin)
				-- ($-1*(Bone)$) node [below left]  {$-\al-\b$};
			\end{scope} 
			\begin{scope}
				\pgftransformcm{1}{0}{-1/2}{sqrt(3)/2}{\pgfpoint{0cm}{0cm}} 
				\draw[style=help lines,dashed] (-6,-6) grid[step=2cm] (6,6);
			\end{scope}
		\end{tikzpicture}
		\caption{Borel-type submonoid of $A_2$ \label{fig1}}
	\end{figure}
\end{example}

\subsection{Borel-type submonoids and the root system}\label{Sec:Borelandrootsystem}

The following subset of $L$ is called the set of {\bf roots} in the even lattice $L$ \cite{CS88}:
\begin{equation}
	\Phi(L)=\{\al\in L: (\al|\al)=2 \}.
\end{equation}
Note that $\Phi(L)$ could be empty. For instance, when $L=\Lambda$ is the Leech lattice \cite{FLM}. 

It is well-known that $\Phi(L)$ is a root system of ADE-type if it is nonempty \cite{CS88}. In this subsection, we justify our term ``Borel-type submonoid'' by relating it to the root system $\Phi(L)$. 

\begin{prop}\label{prop:rootBorel1}
	Let $B\leq L$ be a Borel-type submonoid, and assume that $\Phi(L)\neq \emptyset$. Then $\Phi^+:=B\cap \Phi(L)$ is a set of positive roots in the root system $\Phi(L)$.  
\end{prop}
\begin{proof}
	Recall that a subset $\Phi^+$ of a root system $ \Phi$ is a set of {\em positive roots} if the following conditions hold, see \cite{CS88,FH91,Hum2}: 
	\begin{enumerate}
		\item For all $\al\in \Phi$, exactly one of $\al$ and $-\al$ is in $\Phi^+$.
		\item For $\al,\b\in \Phi^+$ such that $\al+\b\in \Phi$, one has $\al+\b\in \Phi^+$. 
	\end{enumerate}
	Since the root system satisfies $\Phi(L)=-\Phi(L)$, we have $\Phi(L)\cap(-B)=-(\Phi(L)\cap B)=-\Phi^+$. Then by Definition~\ref{def:submonoidsofaPDElattice}, we have $\Phi(L)=(\Phi(L)\cap B)\cup (\Phi(L)\cap (-B))=\Phi^+\cup (-\Phi^+)$, and $\Phi^+\cap (-\Phi^+)\ssq B\cap (-B)=0$. In particular, for any $\al\in \Phi(L)$, either $\al\in \Phi^+$ or $-\al\in \Phi^+$. 
	
	Moreover, let $\al,\b\in \Phi^+$ such that $\al+\b\in \Phi(L)$. Since $\al,\b\in B$, and $B\leq L$ is closed under addition, we have $\al+\b\in B\cap \Phi(L)=\Phi^+$. This shows $\Phi^+=\Phi(L)\cap B$ is a set of positive roots in $\Phi(L)$.  
\end{proof}

There is a converse of Proposition~\ref{prop:rootBorel1} for certain even lattices $L$. In other words, a set of positive roots in $\Phi(L)$ is also contained in some Borel-type submonoid of $L$. First, we let $Q$ be a root lattice of type $A,D$, or $E$. Then the set of roots $\Phi(Q)=\{\al\in Q: (\al|\al)=2\}$ is the root system of type $A,D$, or $E$, respectively, and $Q=\Z\Phi$. See \cite{Hum2,FH91}.  

\begin{lm}\label{lm:rootBorel2}
	Let $Q$ be a root lattice of type $A,D$, or $E$, and let $\Phi^+$ be a set of positive roots in $\Phi(Q)$. Then there exists a Borel-type submonoid $B\leq Q$ such that $\Phi^+ =B\cap \Phi(L)\subset B$. 
\end{lm}
\begin{proof}
	Let $V$ be the Euclidean space spanned by $Q$. By \cite[Theorem 10.1]{Hum2}, there exists a hyperplane $P(\ga)=\{v\in V: (v|\ga)=0\}\subset V$ such that $\Phi^+=\Phi^+(\ga)=\{ \al\in \Phi: (\al|\ga)>0\}$. i.e., $\Phi^+=P^+(\ga)\cap \Phi(L)\subset P^+(\ga)\cap Q$. If $P(\ga)\cap Q=\{0\}$, then $B:=( P^+(\ga)\cap Q)\cup\{0\}$ is a Borel-type submonoid containing $\Phi^+$ such that $B\cap \Phi(L)=P^+(\ga)\cap \Phi(L)=\Phi^+$, we are done. 
	
	Now assume $P(\ga)\cap Q\neq \{0\}$. Since $P(\ga)\cap Q$ is a abelian subgroup of the free abelian group $Q$, it is necessarily free. Hence $P(\ga)\cap Q$ is a lattice in the Euclidean space $W=\spn_\R \{P(\ga)\cap Q\}$ which is contained in $ P(\ga)$. By Proposition~\ref{prop:existenceofBorel}, the lattice $P(\ga)\cap Q$ has a Borel-type submonoid $B_1$. Let $B:=(P^+(\ga)\cap Q)\sqcup B_1$. Since $Q\cap P(\ga)=B_1\cup (-B_1)$, then we have  
	\begin{align*}
		Q&=(P^+(\ga)\cap Q)\sqcup (Q\cap P(\ga))\sqcup (P^-(\ga)\cap Q)\\
		&=\left((P^+(\ga)\cap Q)\sqcup B_1\right)\cup \left((-B_1)\sqcup (P^-(\ga)\cap Q)\right)\\
		&=B\cup (-B). 
	\end{align*}
	Moreover, $B\cap (-B)=B_1\cap (-B_1)=\{0\}$ by the definition of $B_1$, and $B\subset P^{\geq 0}(\ga)$ since $B_1\subset P(\ga)$. Hence $B$ is a Borel-type submonoid of $Q$. Moreover, since $B_1\cap \Phi(L)\ssq P(\ga)\cap \Phi(L)=\emptyset$, we have $B\cap \Phi(L)=(P^+(\ga)\cap \Phi(L))\sqcup (B_1\cap \Phi(L))=\Phi^+\sqcup \emptyset =\Phi^+$. 
\end{proof}

Now let $L$ be an arbitrary positive-definite even lattice in $E$, and assume that $\Phi(L)\neq \emptyset$. Let $Q(L)=\Z \Phi(L)=\Z\al_1\op\ds\op \Z\al_k$ be the root lattice spanned by the root system $\Phi(L)$, where $\Delta(L)=\{\al_1,\ds,\al_k\}$ is a basis of the root system $\Phi(L)$.  

\begin{prop}\label{prop:rootBorel3}
	With the notations as above, assume $\Delta(L)$ can be extended to a basis of $L$. i.e., $L=Q(L)\op \Z \al_{k+1}\op\ds \op \Z\al_{r}$, where $\Delta(L)\cup \{\al_{k+1},\ds,\al_r\}$ is a $\R$-basis of the Euclidean space $E$. Let $\Phi^+$ be a set of positive roots in $\Phi(L)$. Then there exists a Borel-type submonoid $B\leq L$ such that $\Phi^+= B\cap \Phi(L)\subset B$.
\end{prop}
\begin{proof}
	We use induction on $r-k$ to show the existence of $B$. If $L=Q(L)$, the conclusion follows from Lemma~\ref{lm:rootBorel2}. Consider the case when $L=Q(L)\op \Z\al_{k+1}$. Let $V:=\spn_\R Q(L)\subset E$, then $\mathrm{codim}_EV=1$. Since $(\cdot|\cdot)$ is positive-definite, we may choose a vector $\ga\in E\bs\{0\}$ such that $V\perp \ga$. It follows that $V\ssq P(\ga)$, and so $V=P(\ga)$ for the dimension reason. 
	
	We claim that $L\cap V=Q(L)$. Indeed, clearly $Q(L)\ssq L\cap V$. Conversely, let $\al=\sum_{i=1}^k m_i\al_i+m_{k+1}\al_{k+1}=\sum_{j=1}^k\la_j\al_j\in L\cap V$, where $m_i\in \Z$ and $\la_j\in \R$ for all $i,j$. Since $\Delta(L)\cup \{\al_{k+1}\}$ are $\R$-linearly independent, we have $m_{k+1}=0$ and $\al=\sum_{i=1}^k m_i\al_i\in Q(L)$. 
	
	Since $\Phi(L)$ is a root system of $A,D$, or $E$ type, then by Lemma~\ref{lm:rootBorel2}, there exists a Borel-type submonoid $B$ of the lattice $Q(L)$ such that $\Phi^+=B\cap \Phi(L)$. Now let $\tilde{B}:=B\sqcup (L\cap P^+(\ga))$ which is a submonoid in $L$. We have $\tilde{B}\subset P^{\geq 0}(\ga)$ since $B\subset V=P(\ga)$. Moreover,
	\begin{align*}
		L=L\cap E&=(L\cap P^+(\ga))\sqcup Q(L)\sqcup (L\cap -P^+(\ga))\\
		&= \left((L\cap P^+(\ga))\sqcup B\right) \cup \left( (-B)\sqcup (L\cap -P^+(\ga))\right)\\
		&=\tilde{B}\cup (-\tilde{B}),
	\end{align*}  
	and $\tilde{B}\cap (-\tilde{B})=B\cap (-B)=0$. Hence $\tilde{B}\leq L$ is a Borel-type submonoid. Note that $(L\cap P^+(\ga))\cap \Phi(L)=0$ since $\Phi(L)\subset Q(L)\subset  V=P(\ga)$. It follows that $\tilde{B}\cap \Phi(L)=B\cap \Phi(L)=\Phi^+$. This finishes the base case. 
	
	Suppose the conclusion holds for smaller $r-k$. Consider the sublattice $L'=Q(L)\op \Z \al_{k+1}\op\ds \op \Z\al_{r-1}$ in the Euclidean subspace $V=\spn_\R L'$ of $E$. Since $\Delta(L)\cup \{\al_{k+1},\ds,\al_r\}$ is a $\R$-basis $E$, we have $\mathrm{codim}_EV=1$. Similarly, we can find a vector $\ga\in E\bs\{0\}$ such that $V=P(\ga)$ and show that  $L\cap V=L'$. By the induction hypothesis, there exists a Borel-type submonoid $B\leq L'$ such that $\Phi^+=B\cap \Phi(L)$. Similar to the discussion above, we can show that $\tilde{B}:=B\cup P^+(\ga)$ is a Borel-type submonoid of $L$ such that $\tilde{B}\cap \Phi(L)=B\cap \Phi(L)=\Phi^+$ since $\Phi(L)\subset L'\subset P(\ga)$. 
\end{proof}


\section{Borel-type and parabolic-type subVOAs of $V_L$}\label{Sec:3}

In this Section, we first review the construction of lattice VOAs in \cite{FLM} and some related results, then introduce the notions of conic-type, Borel-type and parabolic-type subVOAs of a lattice VOA $V_L$ based on the submonoids of $L$ we introduced in Section~\ref{Sec2}. 

\subsection{Basics of the lattice VOAs}\label{sec2.1}

For the general definitions of vertex operator algebras (VOAs), modules over VOAs and rationality, we refer to \cite{FLM,FHL,DL,LL,Z}. 

Let $L$ be a positive definite even lattice of rank $d\geq 1$, equipped with $\Z$-bilinear form $(\cdot|\cdot):L\times L\ra \Z$. Let $\h:=\C\otimes_{\Z} L$, extend $(\cdot|\cdot)$ to a nondegenerate $\C$-bilinear form $(\cdot|\cdot):\h\times \h\ra \C$, and let $M_{\hat{\h}}(1,0)$ be the level-one Heisenberg VOA associated to $\hat{\h}=\h\otimes \C[t,t^{-1}]\op \C K$. 

We give a brief recap of the construction of the Heisenberg VOA $M_{\hat{\h}}(1,0)$ and its irreducible modules $M_{\hat{\h}}(1,\la)$, which will be used later. The Lie bracket on the affine algebra $\hat{\h}$ is given by:
\begin{equation}\label{2.3'}
	[h_1(m),h_2(n)]=m\delta_{m+n,0} (h_1|h_2) K,\quad \forall h_1,h_2\in \h,\ \text{and}\  m,n\in\Z,
\end{equation}
where we denote $h\o t^m\in \hat{\h}$ by $h(m)$. Then $\hat{\h}=(\hat{\h})_{+}\op (\hat{\h})_{0}\op (\hat{\h})_{-}$, where $(\hat{\h})_{\pm}=\bigoplus_{n\in \Z_{\pm}} \h\o \C t^{n} $, and $(\hat{\h})_{0}=\h\o \C 1\op \C K$. Let $(\hat{\h})_{\geq 0}=(\hat{\h})_{+}\op (\hat{\h})_{0}$, which is a Lie sub-algebra of $\hat{\h}$. 

For each $\la\in \h$, let $e^\la$ be a formal symbol associated to $\la$. Then $\C e^\la$ is a module over $(\hat{\h})_{\geq 0}$, with the module actions given by $h(0)e^\la=(h|\la)e^\la$, $K.e^\la=e^\la$, and $h(n)e^\la=0$, for all $h\in \h$ and $n>0$. Then define $ M_{\hat{\h}}(1,\la)$ to be the induced module:
\begin{equation}\label{2.4'}
	M_{\hat{\h}}(1,\la):=\Ind_{\hat{\h}_{\geq 0}}^{\hat{\h}} \C e^\la=U(\hat{\h})\o _{U(\hat{\h}_{\geq 0})} \C e^\la.
\end{equation}
Then we have $ M_{\hat{\h}}(1,\la)\cong U(\hat{\h}_{<0})\o_{\C} \C e^\la=\spn\{ h_1(-n_1)\ds h_k(-n_k)e^\la: k\geq 0,\ h_1,\ds,h_k\in \h,\ n_{1}\geq \ds \geq n_k\geq 0 \}$ as vector spaces. It was proved in \cite{FLM} that $M_{\hat{\h}}(1,0)$ is a VOA called the level-one Heisenberg VOA, and $M_{\hat{\h}}(1,\la)$, with different $\la\in \h$, are all the irreducible modules over $M_{\hat{\h}}(1,0)$ up to isomorphism.

Let $\epsilon: L\times L\ra \<\pm 1\>$ be a $2$-cocycle of the abelian group $L$ such that $\epsilon(\al,\b)\epsilon(\b,\al)=(-1)^{(\al|\b)}$, for any $\al,\b\in L$. Write $\C^\ep[L]=\bigoplus_{\al\in L} \C e^{\al}$, where $e^{\al}$ is a formal symbol associated to $\al$ for each $\al\in L$ ($e^\al$ was denoted by $\iota(\al)$ in \cite{FLM}). Let $V_L=M_{\hat{\h}}(1,0)\otimes_{\C} \C^\epsilon[L]$, then by the discussion above, we have the following spanning set of $V_L$: 
$$V_L=\spn\{ h_1(-n_1)\ds h_k(-n_k)e^\al: k\geq 0, \al\in L, \ h_1,\ds,h_k\in \h,\ n_{1}\geq \ds \geq n_k\geq 0 \},$$
where we omit the tensor sign $\o$ in the term $h_1(-n_1)\ds h_k(-n_k)\o e^\al$. The vertex operator $Y:V_L\ra \End(V_L)[[z,z^{-1}]]$ has the following definition on the spanning elements of $V_L$: 
\begin{align}
	&Y(h(-1)\vac,z):=h(z)=\sum_{n\in \Z} h(n) z^{-n-1}\quad \left(h(n)e^\al:=0,\ \forall n>0,\ h(0)e^\la:=(h|\al) e^\al\right), \label{2.3}\\
	&Y(e^{\al},z):=E^{-}(-\al, z)E^{+}(-\al,z)e_{\al} z^{\al} \quad \left(z^\al (e^\b):=z^{(\al|\b)}e^\b,\ e_\al(e^\b):=\epsilon (\al,\b)e^{\al+\b}\right),\label{2.4}\\
	&Y(h_1(-n_{1}-1)\dots h_k(-n_{k}-1)e^{\al},z):={\tiny\begin{matrix}\circ \\\circ\end{matrix}}(\partial_{z}^{(n_1)}h_1(z))\dots (\partial_{z}^{(n_k)}h_k(z))Y(e^{\al},z){\tiny\begin{matrix}\circ \\\circ\end{matrix}},\label{2.5}
\end{align}
for any $k\geq 1$, $n_{1}\geq \ds\geq n_k\geq 0$, $h,h_1,\dots h_k\in \h$, and $\al,\b\in L$, where the operators $E^{\pm}$, $\partial^{(n)}_{z}$, $e^\al$, and $z^\al$ in \eqref{2.3}--\eqref{2.5} are given as follows:  
$$ E^{\pm}(-\al,z)=\exp\left(\sum_{n\in \Z_{\pm}}\frac{-\al(n)}{n}z^{-n}\right),\quad \partial_z^{(n)}=\frac{1}{n!}\frac{\partial^n}{\partial z^n},\quad e_\al(e^{\b})=\epsilon(\al,\b)e^{\al+\b},\quad z^\al(e^\b)=z^{(\al|\b)}z^{\b}.
$$
The normal ordering in \eqref{2.5} rearranges the terms in such a way that the right hand side of \eqref{2.5} has the following expression: 
\begin{equation}\label{2.6}
	\sum_{m_1>0,\ds, m_k>0}\sum_{n_1\geq 0,\ds, n_k\geq 0} c_{m_1,\ds,n_k} h_1(-m_1)\ds h_k(-m_k)E^{-}(-\al,z)e_\al z^\al E^{+}(-\al,z)h_1(n_1)\ds h_k(n_k).
\end{equation}
Let $\{\al_1,\dots,\al_d\}$ be an orthonormal basis of $\h$, and let $\om =\frac{1}{2}\sum_{i=1}^d \al_i(-1)^{2}\vac \in M_{\hat{\h}}(1,0)\subset  V_L$.

It was proved in \cite[Appendix A.2]{FLM} that $(V_L,Y,\om ,\vac)$ is a VOA, and $M_{\hat{\h}}(1,0)\leq V_L$ a subVOA with the same Virasoro element $\om$ and vacuum element $\vac$. $V_L$ can be decomposed as a direct sum of irreducible modules over the Heisenberg subVOA $M_{\hat{\h}}(1,0)$ \cite{FLM,D}: 
\begin{equation}\label{2.7}
	V_{L}=\bigoplus_{\al\in L} M_{\hat{\h}}(1,\al),
\end{equation}
where $M_{\hat{\h}}(1,\al)=M_{\hat{\h}}(1,0)\otimes \C e^\al$ for each $\al\in L$.



We shall need the following notion introduced by Huang and Lepowsky in \cite{HL}:

\begin{df}\label{df2.1}
	A \name{vertex algebra without vacuum} is a triple $(V,Y,D)$, where $V$ is a vector spaces, $Y:V\ra \End(V)[[z,z^{-1}]],\ a\mapsto Y(a,z)=\sum_{n\in \Z} a_nz^{-n-1}$ is a linear map, and $D:V\ra V$ is a linear map, satisfying the following axioms:
	\begin{enumerate}
		\item (Truncation property) For any $a,b\in V$, $a_nb=0$ for $n\gg 0$. 
		\item ($D$-derivative property) For any $a\in V$, $[D,Y(a,z)]=\frac{d}{dz}Y(a,z)$.
		\item (Skew-symmetry) For any $a,b\in V$, $Y(a,z)b=e^{zD}Y(b,-z)a$. 
		\item (Jacobi identity) For any $a,b,c\in V$, 
		\begin{align*}& \displaystyle{z^{-1}_0\delta\left(\frac{z_1-z_2}{z_0}\right)
				Y(a,z_1)Y(b,z_2)-z^{-1}_0\delta\left(\frac{-z_2+z_1}{z_0}\right)
				Y(b,z_2)Y(a,z_1)}
			\displaystyle{=z_2^{-1}\delta
				\left(\frac{z_1-z_0}{z_2}\right)
				Y(Y(a,z_0)b,z_2)}.
		\end{align*}
	\end{enumerate}
	In particular, a VOA $(V,Y,\vac,\om)$ is a vertex algebra without vacuum with $D=L(-1)$. A sub-vertex algebra without vacuum $W\leq V$ is a subspace $W\subset V$ that is closed under $Y$ and $L(-1)$. 
\end{df}

\subsection{SubVOAs of $V_L$ associated to submonoids of $L$}
First, we observe that submonoids and semigroups of $L$ give rise to sub-vertex algebras (without vacuum) of $V_L$: 
\begin{prop}\label{prop2.3}
	Let $L$ be a positive-definite even lattice, $M\leq L$ be a submonoid, and $S\subset L$ be a sub-semigroup. Then 
	\begin{enumerate}
		\item $(V_M:=\bigoplus_{\al\in M} M_{\hat{\h}}(1,\al),Y,\om ,\vac)$ is a CFT-type subVOA of $(V_L,Y,\om ,\vac)$. 
		
		\item $(V_S:=\bigoplus_{\al\in S} M_{\hat{\h}}(1,\al),Y,L(-1))$ is a sub-vertex algebra without vacuum of $(V_L,Y,L(-1))$. If, furthermore, $S\subset M$ and $M+S\ssq S$, then $V_S$ is an ideal of $V_M$. 
	\end{enumerate}
	We call $(V_M,Y,\om ,\vac)$ (resp. $(V_S,Y,L(-1))$) the \name{subVOA (resp. sub-vertex algebra without vacuum) of $V_L$ associated to $M$ (resp. $S$)}. 
\end{prop}
\begin{proof}
	For any $\al,\b\in M$, 	by \eqref{2.3} and \eqref{2.4},	we have 
	\begin{align*}
		Y(e^{\al},z)e^{\b}&=E^{-}(-\al, z)E^{+}(-\al,z)e_{\al} z^{\al}(e^{\b})=E^{-}(-\al, z)E^{+}(-\al,z)\epsilon(\al,\b) z^{(\al|\b)} e^{\al+\b},\\
		&=\exp(\sum_{n<0}-\frac{\al(n)}{n} z^{-n})\epsilon(\al,\b) z^{(\al|\b)} e^{\al+\b}.
	\end{align*}
	The right-hand-side is contained in $M_{\hat{\h}}(1,\al+\b)((z))\subset V_M((z))$, in view of the decomposition \eqref{2.7}. More generally, for any $h_1(-n_{1}-1)\dots h_k(-n_{k}-1)e^{\al}\in M_{\hat{\h}}(1,\al)$ and $h'_1(-m_{1}-1)\dots h'_r(-m_{r}-1)e^{\b}\in M_{\hat{\h}}(1,\b)$, with $\al,\b\in M$, it is easy to see from \eqref{2.5} and \eqref{2.6} that
	\begin{align*}
		Y(h_1(-n_{1}-1)\dots h_k(-n_{k}-1)e^{\al},z)h'_1(-m_{1}-1)\dots h'_r(-m_{r}-1)e^{\b}\in M_{\hat{\h}}(1,\al+\b)((z)).
	\end{align*}
	Since $M$ is closed under addition and $M_{\hat{\h}}(1,0)\subset  V_M$, it follows that $V_M$ is a sub-VOA of $V_L$. Since $V_M$ has the same Virasoro element as $V_L$, we have $(V_M)_n\ssq (V_L)_n$ for each $n\geq 0$, and $(V_M)_0=(V_L)_0=\C\vac$. Thus $V_M$ is of CFT-type. The second statement is also clear since $S$ is closed under addition, and $L(-1) M_{\hat{\h}}(1,\al)\ssq M_{\hat{\h}}(1,\al)$ for any $\al\in S$. 
\end{proof}

\begin{remark}
	The proof of Proposition~\ref{prop2.3} essentially depends on the fact that 
	\begin{equation}\label{2.8}
		Y(M_{\hat{\h}}(1,\al),z)M_{\hat{\h}}(1,\b)\subset M_{\hat{\h}}(1,\al+\b)((z)),\quad \al,\b\in L,
	\end{equation}
	where $Y$ is the vertex operator of the lattice VOA $V_L$. In other words, modules over the Heisenberg VOA are all simple current. When $L=\Z\al$, it was observed by Dong that $V_{\N \al}$ is a subVOA of $V_{\Z\al}$, see \cite[Proposition 4.1]{D}. Proposition \ref{prop2.3} is a higher rank generalization of this result, noting that $\N\al$ is a submonoid of $\Z\al$. 
\end{remark}



Using the different types of submonoids in $L$ introduced in Definition~\ref{def:submonoidsofaPDElattice}, together with Proposition~\ref{prop2.3}, we give the following definition of subVOAs in $V_L$:

\begin{df}\label{def:subVOAsofVL}
Let $L$ be a positive-definite even lattice, and let $C\leq L$, $B\leq L$, and $P\leq L$ be a conic-type, Borel-type, and parabolic-type submonoid of $L$, respectively, as in Definition~\ref{def:submonoidsofaPDElattice}. Then 
\begin{enumerate}
	\item $V_C=\bigoplus_{\al\in C} M_{\hat{\h}}(1,\al)$ is called a \name{conic-type subVOA} of $V_L$. 
	\item  $V_B=\bigoplus_{\al\in B}M_{\hat{\h}}(1,\al)$ is called a \name{Borel-type subVOA} of $V_L$. 
	\item $V_P=\bigoplus_{\al\in P}M_{\hat{\h}}(1,\al)$ is called a \name{parabolic-type subVOA} of $V_L$. 
\end{enumerate}
\end{df}

Let $V_B=\bigoplus_{\al\in B}M_{\hat{\h}}(1,\al)$ be a Borel-type subVOA. 
We may view $M_{\hat{\h}}(1,0)\leq V_B$ as the analog of ``Cartan part'', and view $M_{\hat{\h}}(1,\al)$ with $\al\in B\bs\{0\}$ as the analog of a ``root space'' associated to a ``positive root'' $\al\in B$. However, unlike the Lie algebra case, the ``Cartan subalgebra'' $M_{\hat{\h}}(1,0)$ is not commutative and the ``root space'' $M_{\hat{\h}}(1,\al)$ is not one-dimensional. 

\begin{remark}We remark some easy facts from our definition: 
\begin{enumerate}
	\item The subVOAs $V_C,V_B$ and $V_P$ are of CFT-type. They also have the same vacuum $\vac$ and Virasoro element $\om$ with the lattice VOA $V_L$, see  Proposition~\ref{prop2.3}.
	\item Any Borel-type subVOA $V_B$ contains a conic-type subVOA $V_C$. Any Borel-type subVOA is also of parabolic-type, see Proposition~\ref{prop:propertiesofBorel}. 
\end{enumerate}
\end{remark}

\begin{example}	 \label{ex2.7}
We have natural examples of conic, Borel, and parabolic-type subVOAs of low-rank lattices arising from Examples~\ref{ex:rankoneBorel} and \ref{ex:ranktwoparabolic}. 
\begin{enumerate}
	\item 	$V_B=V_{\Z_{\geq 0}\al}$ is a conic and Borel-type subVOA of $V_{\Z\al}$. Note that $V_{\Z\al}=V_{B}\op V_{\Z_{<0}\al}$, where  $ V_{\Z_{<0}\al}$ is a sub-vertex algebra without vacuum of $V_{\Z\al}$. 
	
	\item Let $L=A_2$, and let $B$ be given by \eqref{eq:exampleofB}. Then $V_B$ is a Borel-type subVOA of $V_{A_2}$. Moreover, $V_{\Z\al\op \Z_{\geq 0}\b}$ is a parabolic-type subVOA of $V_{A_2}$ that is not of Borel-type. $V_P$ gives rise to a decomposition $V_{A_2}=V_{P}\op V_{\Z\al\op \Z_{<0}\b}$, where $V_{\Z\al\op \Z_{<0}\b}$ is a sub-vertex algebra without vacuum of $V_{A_2}$. 
\end{enumerate}
\end{example}


\begin{remark}	The decompositions $V_{\Z\al}=V_{B} \op V_{\Z_{< 0} \al}$ and $V_{A_2}=V_{P}\op V_{\Z\al\op \Z_{<0}\b}$ in Example~\ref{ex2.7} first appeared in the study of Rota-Baxter operators and classical Yang-Baxter equations for VOAs by Bai, Guo, the author, and Wang in \cite{BGL,BGLW}. The projections $p:V_{\Z\al}\ra V_B\subset V_{\Z\al}$ and $p:V_{A_2}\ra V_P\subset V_{A_2}$ are natural examples of weight $-1$ Rota-Baxter operator for VOAs, which can be viewed as the operator-form classical Yang-Baxter equation on VOAs.
\end{remark}


\section{Basic properties of conic, Borel, and parabolic-type subVOAs} \label{Sec:4}

In this Section, we discuss the basic properties of the conic, Borel, parabolic-type subVOAs. We prove that these CFT-type VOAs are irrational and some of them are $C_1$-cofinite. We also show that they share many similar properties as the usual Borel and parabolic Lie subalgebras of a semisimple Lie algebra, which further justifies our terminologies. 

\subsection{Relations with the root system and degree-one Lie algebra of $V_L$}\label{Sec:4.1}
In this subsection, we justify our terminologies of Borel and parabolic-type subVOAs in Definition~\ref{def:subVOAsofVL} by relating them with the root system $\Phi(L)$ of the lattice $L$ and the degree-one Lie algebra of a lattice VOA. 

For any CFT-type VOA $V=\bigoplus_{n=0}^\infty V_n$, recall that the degree-one subspace $\g=V_1$ is a Lie algebra with respect to the Lie bracket: 
\begin{equation}
[a,b]=a_0b=\Res_z Y(a,z)b,\quad a,b\in V_1,
\end{equation}
see  \cite{B,FLM}. It was proved by Dong and Mason that $\g=V_1$ is a reductive Lie algebra if the VOA $V$ is {\bf strongly rational} i.e., if $V$ is  of CFT-type, rational, $C_2$-cofinite, $L(1)V_1=0$, and $V\cong V'$ as $V$-modules \cite{DM04}. 

In particular, let $V_L$ be a lattice VOA, then $\g=(V_L)_1$ is a reductive Lie algebra since $V_L$ is strongly rational \cite{D}. Recall that $e^\al\in (V_L)_1$ if and only if $(\al|\al)=2$. i.e., $\al\in \Phi(L)$. Then $\g=\h$ is abelian if $\Phi(L)=\emptyset$. On the other hand, if $\Phi(L)\neq \emptyset$, then 
\begin{equation}\label{eq:thefirstlevelLiealgebra}
\g=\mathfrak{s}\op \mathfrak{c}=\h+ \sum_{\al \in \Phi(L)} \C e^\al,
\end{equation}
where $\mathfrak{s}=H+\sum_{\al\in \Phi(L)}\C e^{\al}$ is a semisimple Lie subalgebra of $\g$ associated to the ADE-type root system $\Phi(L)$, with the Cartan subalgebra $H=\spn\{\al(-1)\vac: \al\in \Phi(L)\}$, $\C e^\al= \g_{\al}$ is the root space associated to the root $\al\in \Phi(L)$, $\mathfrak{c}$ is the center of $\g$, and $\h=H+\mathfrak{c}$. 

Using the relations between Borel-type submonoids $B\leq L$ and the root system $\Phi(L)$ in Section~\ref{Sec:Borelandrootsystem}, we have can prove the following facts: 

\begin{prop}\label{prop:Borelfirstlevel1}
Let $V_P$ (resp. $V_B$) be a parabolic-type (resp. Borel-type) subVOA of $V_L$, and assume that $\Phi(L)\neq \emptyset$. Then $\fp=(V_P)_1$ (resp. $\fb=(V_B)_1$) is a parabolic (resp. Borel) subalgebra of the reductive Lie algebra $\g=(V_L)_1$.  
\end{prop}
\begin{proof}
Let $P\leq L$ be a parabolic-type submonoid containing a Borel-type submonoid $B$. By Proposition~\ref{prop:rootBorel1}, $B\cap \Phi(L)=\Phi^+$ is a set of positive roots in $\Phi(L)$. Since $\wt (e^\al)=(\al|\al)/2>1$ for any $\al\in B\bs \Phi(L)$, we have 
\begin{equation}\label{eq:firstlevelofBorel}
	(V_B)_1=\bigoplus_{\al\in B} \left(M_{\hat{\h}}(1,\al)\right)_1=\spn\{\al(-1)\vac: \al\in \h\}+ \sum_{\al\in B\cap \Phi(L)} \C e^{\al}=\h+\sum_{\al\in \Phi^+}\C e^\al,
\end{equation}
which is a Borel Lie subalgebra $\fb$ of the reductive Lie algebra $\g$ \cite{B1,FH91}. Moreover, since $B\ssq P$, we have $\Phi^+\ssq P\cap \Phi(L)$. Then $P\cap \Phi(L)$ is a subset of the root system $\Phi(L)$ containing all the positive roots and a subset of negative simple roots. Hence
\begin{equation}\label{eq:firstlevelofparabolic}
	(V_P)_1=\bigoplus_{\al\in P} \left(M_{\hat{\h}}(1,\al)\right)_1=\h+\sum_{\al\in P\cap \Phi(L)} \C e^\al
\end{equation}
is a parabolic Lie subalgebra $\fp$ of $\g$ \cite{B1,BT,FH91}. 
\end{proof}

There is also a converse of Proposition~\ref{prop:Borelfirstlevel1} that follows from Lemma~\ref{lm:rootBorel2} and Proposition~\ref{prop:rootBorel3}. 

\begin{prop}\label{prop:Borelfirstlevel2}
Assume that $\Phi(L)\neq \emptyset$, and let $\Delta(L)$ be a basis of the root system $\Phi(L)$. Assume $\Delta(L)$ can be extended to a basis of $L$ as in  Proposition~\ref{prop:rootBorel3}.

Let $\fb$ (resp. $\fp$) be a Borel (resp. parabolic) Lie subalgebra of the reductive Lie algebra $\g=(V_L)_1$ that contains $\h=\spn\{h(-1)\vac: h\in \h\}\subset (V_L)_1$. Then there exists a Borel-type subVOA $V_B\leq V_L$ such that $\fb= (V_B)_1$ and a parabolic-type subVOA $V_P\leq V_L$ such that $\fp\ssq (V_P)_1$. 
\end{prop}
\begin{proof}
Note that a Borel subalgebra $\fb\leq \g$ gives rise to a set of positive roots $\Phi^+$ of $\Phi(L)$ such that $\fb=\h+\sum_{\al\in \Phi^+} \C e^\al$. By Proposition~\ref{prop:rootBorel3}, there exists a Borel-type submonoid $B\leq L$ such that $\Phi^+=B\cap \Phi(L)$. It follows from \eqref{eq:firstlevelofBorel} that $\fb=(V_B)_1$. 

Let $\fp =\h+ \sum_{\al\in \Gamma} \C e^\al$ be a parabolic subalgebra of $\g$ containing $\h$. By definition, $\Gamma=\Phi^+ \cup S$, where $S$ is a subset of $\Phi^-$ generated by some negative simple roots $\{-\al_{i_1},\ds ,-\al_{i_k}\}\subset -\Delta(L)$. Let $B\leq L$ be a Borel-type submonoid such that $\Phi^+=B\cap \Phi(L)$, then
$$P:=B+\Z\al_{i_1}+\ds +\Z\al_{i_k}\subset L$$
is also a submonoid of $L$ that contains $B$, and $P$ is a parabolic-type submonoid. Since $\Phi^+\subset B\cap \Phi(L) \subset P\cap \Phi(L) $ and $S\subset (\Z\al_{i_1}+\ds +\Z\al_{i_k})\cap \Phi(L)\subset P\cap \Phi(L)$, we have $\Gamma=\Phi^+\cup S\subset P\cap \Phi(L)$. It follows from \eqref{eq:firstlevelofparabolic} that 
$\fp=\left(\h+ \sum_{\al\in \Gamma} \C e^\al\right) \ssq \left(\h+\sum_{\al\in P\cap \Phi(L)} \C e^\al\right) =(V_{P})_1.$
\end{proof}


Certain conic-type subVOAs $V_C$ satisfy a conjugacy property that is similar to the conjugacy properties of the Borel-subalgebras of a semisimple Lie algebra $\g$ or the Borel-subgroups of a reductive algebraic group $G$. 


\begin{prop}\label{prop2.10}
Let $Q$ be a root lattice associated to the root system $\Phi$ of type $A,D$, or $E$. Then the conic-type subVOAs of the form $V_{C_\Delta}$, where $\Delta=\{\al_1,\ds,\al_r\}$ is a basis of $\Phi$ and $C_\Delta=\Z_{\geq 0}\al_1\op \ds \op \Z_{\geq 0}\al_r$, are conjugate under $\mathrm{Aut}(V_Q)$. 
\end{prop}
\begin{proof}
Given a lattice isometry $\si\in \mathrm{Aut}(Q)$ \eqref{eq:defofAutL}, we can lift it to an element $\hat{\si}\in \mathrm{Aut}(V_Q)$ as follows: First, extend $\si$ to a $\C$-linear isomorphism $\si:\h\ra \h$ by letting  $\si(\la \al)=\la\si(\al),$ where $\la\in \C$ and $\al\in L$. Then define
\begin{equation}\label{2.10}
	\hat{\si}(h_1(-n_1)\ds h_r(-n_r)e^\al):=(\si h_1)(-n_1)\ds (\si h_r)(-n_r)e^{\si \al},
\end{equation}
where $h_i\in \h$, $n_1\geq \ds\geq n_r\geq 1$, and $\al\in Q$. In particular, $\hat{\si} (M_{\hat{\h}}(1,\al))=M_{\hat{\h}}(1,\si(\al))$ for any $\al\in Q$, see  \cite[Section 2.4]{DN99} for more detail. Let $\Delta$ and $\Delta'$ be two sets of simple roots of $\Phi$. Then there exists $w\in W(\Phi)$, the Weyl group of $\Phi$, such that $w (\Delta)=\Delta'$ \cite{Hum2}. Since $w$ is generated by simple reflections, we have $w\in  \mathrm{Aut}(Q)$ and $w(C_\Delta)=C_{\Delta'}$, see Proposition~\ref{prop:propertiesofBorel}.
Then by \eqref{2.10}, we have $\hat{w}\in \mathrm{Aut}(V_Q)$ and $$\hat{w}(V_{C_{\Delta}})=\hat{w}\left(\bigoplus_{\al\in C_{\Delta}}M_{\hat{\h}}(1,\al)\right)=\bigoplus_{\al\in C_{\Delta}}M_{\hat{\h}}(1,w(\al))=\bigoplus_{\b\in w(C_{\Delta})}M_{\hat{\h}}(1,\b)=V_{C_{\Delta'}}.$$ 
Thus, elements in $\{V_{C_{\Delta}}: \Delta\ \mathrm{is\ a\ basis\ of\ } \Phi\}$ are conjugate under $\mathrm{Aut}(V_Q)$. 
\end{proof}

\subsection{Irrationality and non-simplicity}

In Lie theory, we know that a Borel subalgebra $\fb$ or a proper parabolic-type subalgebra $\fp$ of a semisimple Lie algebra $\g$ is not simple nor semisimple as a Lie algebra \cite{BT,Hum2}. We have similar results for the Borel-type subVOAs.

\begin{prop}\label{prop:nonsimpleofconic}
Conic-type subVOAs of a lattice VOA $V_L$ are {\bf not} simple nor rational. 
\end{prop}
\begin{proof}
Let $C=\Z_{\geq 0} \al_1\op \ds\op \Z_{\geq 0}\al_{r}$ be a conic-type submonoid of $L$, where $\{\al_1,\ds,\al_r\}$ is a basis of $L$. Then $S=\Z_{>0}\al_1\op \ds \op \Z_{>0} \al_r$ is a sub-semigroup of $C$ such that $C+S\ssq S$. By Proposition~\ref{prop2.3}, $V_C$ has a nontrivial ideal $V_S$. Hence $V_C$ is not simple. 

Recall that a rational VOA necessarily has finitely many irreducible modules up to isomorphism \cite{DLM1,Z}. To show $V_C$ is irrational, consider the following sequence of sub-semigroups of $C$: 
\begin{align*}
	C^0:&=C\bs\{0\}=\{m\al_1+m_2\al_2+\ds+m_r\al_r\in C: m\geq 1\ \mathrm{or}\ m=0, (m_2,\ds,m_r)\neq 0\},\\
	C^1:&=C\bs\{0,\al_1\}=\{m\al_1+m_2\al_2+\ds+m_r\al_r\in C:m\geq 2\ \mathrm{or}\ m=0, (m_2,\ds,m_r)\neq 0\},\\
	C^2:&=C\bs\{ 0,\al_1,2\al_1\}=\{m\al_1+m_2\al_2+\ds+m_r\al_r\in C:m\geq 3\ \mathrm{or}\ m=0, (m_2,\ds,m_r)\neq 0\},\\
	&\vdots
\end{align*}
Clearly, $C\supset C^0\supset C^1\supset C^2\supset \ds $. We claim that $C+C^k\ssq C^k$ for all $k\geq 0$. Indeed, suppose there exists $\al=m_1\al_1+m_2\al_2+\ds+m_r\al_r\in C$ and $\b=n_1\al_1+n_2\al_2+\ds+n_r\al_r\in C^k$ such that $\al+\b\notin C^k=C\bs\{0,\al_1,\ds k\al_1\}$. Then  $\al+\b=(m_1+n_1)\al_1+\ds+(m_r+n_r)\al_r=p\al_1$ for some $0\leq p\leq k$. Since $m_i,n_i\geq 0$ for all $i$, we must have $m_1+n_1=p$ and $m_j=n_j=0$ for all $j\geq 2$. But then $\b=n_1\al_1\in \{0,\al_1,\ds k\al_1 \}$, a contradiction. Thus, $C+C^k\ssq C^k$, and $C^k$ is a sub-semigroup of $C$ for all $k$. 

By Proposition~\ref{prop2.3}, $V_{C^0}\supset V_{C^1}\supset V_{C^2}\supset  \ds $ is a sequence of ideals in the VOA $V_C$ such that $$V_C/V_{C^0}\cong M_{\hat{\h}}(1,0),\quad  V_{C^k}/V_{C^{k+1}}\cong  M_{\hat{\h}}(1,(k+1)\al_1),\quad k\geq 0.$$

For non-negative integers $p\neq q$, the $V_C$-modules $M_{\hat{\h}}(1,p\al_1)$ is not isomorphic to $M_{\hat{\h}}(1,q\al_1)$ since they are not isomorphic as $M_{\hat{\h}}(1,0)$-modules. Thus, $V_C$ has infinitely many non-isomorphic irreducible modules. This shows $V_C$ is irrational. 
\end{proof}

By adopting a similar idea as the proof of the previous Proposition, we can prove the non-simplicity and irrationality of the Borel-type and certain parabolic-type subVOAs.

\begin{prop}\label{prop2.7}
Let $P\leq L$ be a parabolic-type submonoid containing $B$ such that $P\subset P^{\geq 0}(\ga)$, where $P(\ga)$ is the unique hyperplane that defines $B$ (cf. Proposition~\ref{prop:propertiesofBorel}). Then $V_P$ is {\bf not} simple nor rational. In particular, Borel-type subVOAs of $V_L$ are {\bf not} simple nor rational. 
\end{prop}
\begin{proof}

Let $S':=P\cap P^+(\ga)=B\cap P^+(\ga)$. It is clear that $S'$ is a nontrivial sub-semigroup of $P$. For any $\al\in S'$ and $\b\in P$, we have $(\al+\b|\ga)=(\al|\ga)+(\b|\ga)>0$ as $(\al|\ga)>0$ and $(\b|\ga)\geq 0$. Hence $\al+\b\in S'$ and so $P+S'\ssq S'$. Also note that $0\in P\bs S'$. Hence $V_P$ has a proper ideal $V_{S'}$ by  Proposition~\ref{prop2.3}. 

Recall that a VOA $V$ is rational if any admissible module $M$ of $V$ is semisimple \cite{DLM1}.
Suppose $V_P$ is rational. Note that $V_{S'}=\bigoplus_{\al\in S'}M_{\hat{\h}}(1,\al)\leq V_P$ is an admissible submodule since each $M_{\hat{\h}}(1,\al)$ with $\al\in S'$ is admissible. Then there exists an admissible submodule $N\leq V_P$ such that $V_P=V_{S'}\op N$. For any $\al\in S'$, since $(\al|\al)/2>0$, we have $V_{S'}\ssq (V_P)_+$. Thus, $\vac\in N$ and $V_P=V_P.\vac\ssq N$. Then $V_{S'}=0$, which is a contradiction. Hence $V_P$ is not rational. 
\end{proof}

\subsection{Normalizer property}
In Lie theory, we know that if $P\leq G$ is a Borel or parabolic subgroup of a linear algebraic group $G$, then $N_{G}(P)=P$; and if $\fb\leq \g$ is a Borel subalgebra of a semisimple Lie algebra $\g$, then $\mathfrak{n}_\g(\fb)=\fb$. See \cite[Chapter 23]{Hum1}. We have a similar result for the parabolic-type subVOAs as well.

Let $(V,Y,\vac,\om)$ be a VOA, and $W\leq V$ be a sub-vertex algebra. There is an analog of the ``centralizer'' of $W$ in $V$ called the commutant defined by Frenkel and Zhu in \cite{FZ}. By definition, $\mathrm{Com}_V(W)=\{ a\in V:  w_ja=0,\ \mathrm{for\ all}\ j\geq 0,\ w\in W\}$.
We define the normalizer of $W$ in $V$ as follows: 

\begin{df}\label{def:normalizer}
Let $(V,Y,\vac,\om)$ be a VOA, and $W\leq V$ be a sub-vertex algebra without vacuum. 
\begin{equation}\label{2.11}
	N_{V}(W):=\{ a\in V: a_j W\ssq W,\ \mathrm{for\ any}\ j\geq 0\}
\end{equation}
is called the \name{normalizer of $W$ in $V$}. 
\end{df}

\begin{lm}\label{lm2.11}
$N_{V}(W)\leq V$ is a sub-vertex algebra of $V$ with $W\ssq N_V(W)$. Moreover, $N_{V}(W)$ can also be characterized as follows:
\begin{equation}\label{2.12}
	N_{V}(W):=\{ a\in V: b_ja \in W,\ \mathrm{for\ any}\ b\in W, j\geq 0\}.
\end{equation}
In particular, we have $\mathrm{Com}_V(W)\ssq N_V(W)$. 
\end{lm}
\begin{proof}
Since $\vac_j W=0$ for any $j\geq 0$, we have $\vac \in N_V(W)$. Since $Y(W,z)W\ssq W[[z,z^{-1}]]$, clearly $W\ssq N_V(W)$. Let $a,b\in N_V(W)$,  $w\in W$, $n\in \Z$, and $m\geq 0$. By the Jacobi identity, we have 
$$(a_nb)_mw=\sum_{j\geq 0}\binom{n}{j} (-1)^j a_{n-j}b_{m+j}w-(-1)^n\sum_{j\geq 0}\binom{n}{j}(-1)^j b_{m+n-j}a_jw\equiv 0\pmod{ W} $$
since $b_{m+j}w\in W$ and $a_jw\in W$ for any $j\geq 0$. Thus, $N_V(W)$ is a sub-vertex algebra of $V$.

Moreover, by the skew-symmetry of $Y$, we have $a_nb=\sum_{i\geq 0} 1/i! (-1)^{n+i+1} L(-1)^ib_{n+i}a$. By the definition of vertex algebra without vacuum \ref{df2.1}, we have $L(-1)W\ssq W$. Let $a\in V$ be a given element. Then $a_{n}w\in W$ for any $n\geq 0$ and $w\in W$ if and only if $w_n a\in W$ for any $n\geq 0$ and $w\in W$. This shows  the definitions of \eqref{2.11} and \eqref{2.12} are equivalent. 
\end{proof}
\begin{prop}\label{prop2.13}
Let $M\leq L$ be an abelian submonoid. Then	$N_{V_L}(V_M)=V_M$. In particular, for any parabolic-type submonoid $P\leq L$, we have $N_{V_L}(V_P)=V_P$. 
\end{prop}
\begin{proof}
Note that $V_L=\bigoplus_{\ga\in L\bs M}M_{\hat{\h}}(1,\ga)\op V_M$, in view of \eqref{2.7}. Let $a=u+v\in N_{V_L}(V_M)$ with $u=\sum_{\ga\in L\bs M} u_\ga$ and $v\in V_M$, where $u_\ga=0$ for all but finitely many $\ga\in L\bs M$. Since $V_M\ssq N_{V_L}(V_M)$ by Lemma~\ref{lm2.11}, we have $u\in N_{V_L}(V_M)$. 

Now we show that $u=0$. Since $0\in M$, we have $M_{\hat{\h}}(1,0)\ssq V_M$. Then by \eqref{2.12}, we have $(h(-1)\vac)_0 u=\sum_{\ga\in L\bs M} h(0)u_\ga=\sum_{\ga\in L\bs M}(h| \ga) u_\ga\in V_M\cap \bigoplus_{\ga\in L\bs M}M_{\hat{\h}}(1,\ga)=0$. Hence $(h|\ga)u_\ga=0$ for all $\ga\in L\bs M$. For a fixed $\ga\in L\bs M$, choose $h\in \h$ s.t. $(h|\ga)\neq 0$, then we have $u_\ga=0$ for all $\ga$, and so $u=\sum_{\ga\in L\bs M}u_\ga=0$. 
\end{proof}

\begin{remark}
In general, the normalizer $N_V(W)$ is not necessarily equal to $W$. For example, let $S\subset L$ be a sub-semigroup and $0\notin S$. By Proposition~\ref{prop2.3}, $V_S\leq V_L$ is a sub-vertex algebra without vacuum. It is easy to see that $M_{\hat{\h}}(1,0)\subset N_{V_L}(V_S)$ but $M_{\hat{\h}}(1,0)\nsubset V_S$. 
\end{remark}

\section{Strongly generation property of conic, Borel, and parabolic-type subVOAs}\label{Sec:5}
In this Section, we discuss the strongly finitely generating property, or equivalently, the $C_1$-cofiniteness of conic, Borel, and parabolic-type subVOAs introduced in Definition~\ref{def:subVOAsofVL}. There is no analog of such property in the classical Lie theory. The reason is that VOAs are structurally modules over certain infinite-dimensional Lie algebras. We show that some of these subVOAs of $V_L$ are $C_1$-cofinite, some of them are not. Hence they give rise to natural examples of CFT-type non-$C_1$-cofinite VOAs, which are rare among the classical examples of CFT-type VOAs.


\subsection{$C_1$-cofiniteness of small rank conic-type subVOAs}\label{Sec:5.1}

The notion of strongly generation property of VOAs was introduced by Kac in \cite{K}. A CFT-type VOA $V$ is called {\bf strongly generated} by a subset $U\ssq V$ if $V$ is spanned by elements of the following form: 
\begin{equation}\label{2.13}
u^{1}_{-n_1}u^2_{-n_2}\dots u^k_{-n_k}u,\quad \mathrm{where}\quad u^1,\dots ,u^k,u\in U,\quad n_{1}\geq n_{2}\geq \dots \geq n_{k}\geq 1.
\end{equation}

Let $V$ be a CFT-type VOA. The subspace $C_1(V)\subset V$ was introduced by Li in \cite{L99}. 
\begin{equation}\label{eq:defofC1}
C_1(V):=\spn\{a_{-1}b:a,b\in V_+ \}+\spn\{L(-1)v: v\in V_+\}.
\end{equation}
$V$ is called $C_1$-cofinite if $\dim V/C_1(V)<\infty$. A similar $C_1$-cofinite condition on $V$-modules was introduced by Huang in \cite{H05}. Karel and Li proved that a VOA $V$ is strongly generated by a finite-dimensional subspace $U\subset V$ if and only if $V$ is $C_1$-cofinite \cite{L1,KL}.

\begin{lm}\label{lm2.14}
Let $C=\Z_{\geq 0} \al_1\op \ds\op \Z_{\geq 0}\al_{r}$ be a conic-type submonoid of $L$. Assume that $(\al_i|\al_j)\geq 0$, for all $1\leq i\neq j\leq r$. Then $V_C$ is strongly generated by $U:=\{\vac,\al_i(-1)\vac, e^{\al_i}:1\leq i\leq r\}.$ In particular, $V_C$ is $C_1$-cofinite.
\end{lm}
\begin{proof}
Let $W$ be the subspace of $V_{\Z_{\geq 0}\al}$ spanned by elements of the form \eqref{2.13},  with $u^j,u\in U$.
We need to show that each $ M_{\hat{\h}}(1,n_1\al_1+\ds+n_r\al_r)$ is contained in $W$, for all non-negative integers $n_j\geq 0$. Clearly, $M_{\hat{\h}}(1,0)\subset W$. Since $M_{\hat{\h}}(1,\al)=M_{\hat{\h}}(1,0)\otimes_\C \C e^{\al}$, and $M_{\hat{\h}}(1,0)$ is strongly generated by $\{\al_{1}(-1)\vac,\ds,\al_{r}(-1)\vac\}$, we only need to show $e^{n_1\al_1+\ds+n_r\al_r}\in W$. 

The following claim is very useful for our discussion in this Section: 
\begin{equation}\label{2.14}
	\mathrm{If}\ e^\al\in W\ \mathrm{and}\ e^\b\in W\ \mathrm{such\ that\ }(\al|\b)\geq 0,\ \ \mathrm{then}\ e^{\al+\b}\in W.  
\end{equation}
To prove \eqref{2.14}, we observe that 
\begin{align*}
	e^{\al}_{-(\al|\b)-1} e^{\b} &=\Res_{z} z^{-(\al|\b)-1} E^{-}(-\al,z)E^{+}(-\al,z)e_{\al} z^{\al} e^{\b}\\
	&=\Res_{z} z^{-(\al|\b)-1} E^{-}(-\al,z) z^{(\al|\b)} \epsilon(\al,\b) e^{\al+\b}\\
	&=\epsilon(\al,\b) e^{\al+\b}\equiv 0\pmod{W}.
\end{align*}
Hence $e^{\al+\b}\in W$ as $\epsilon(\al,\b)\neq 0$. Furthermore, since $(\al_i|\al_j)\geq 0$ for all $1\leq i\neq j\leq r$ and $(\al_i|\al_i)=2N_i>0$ for all $i$,  we have
$(m_1\al_1+\ds+m_r\al_r|n_1\al_1+\ds+n_r\al_r)=\sum_{i,j=1}^r m_in_j(\al_i|\al_j)\geq 0,$
for any $m_i,n_j\geq 0$. In particular, if $e^{m_1\al_1+\ds+m_r\al_r}\in W$ and $e^{n_1\al_1+\ds+n_r\al_r}\in W$, we have $e^{(n_1+m_1)\al_1+\ds+(n_r+m_r)\al_r}\in W$ by \eqref{2.14}. Then it follows from an easy induction that $e^{n_1\al_1+\ds+n_r\al_r}\in W$ for any $n_i\geq 0$.  
\end{proof}

It follows immediately from the Lemma and Example~\ref{ex:rankoneBorel} that any rank-one conic-type or Borel-type VOA $V_{B}$ is $C_1$-cofinite. However, the rank-two case is much more complicated. We have the following general theorem:

\begin{thm}\label{thm:ranktwocoincC1}
Let  $L=\Z\al\op \Z\b$ be a rank-two positive-definite even lattice such that $(\al|\b)=-n$, $(\b|\b)=2k$ and $(\al|\al)=2\ell$, with $n\geq 1$, and $k\geq \ell\geq 1$. Suppose the integers $n,k,\ell$ satisfy following condition:
\begin{equation}\label{eq:conditionforranktwoC_1}
	n^2+\ell^2-4\ell k\leq 0.
\end{equation}
Then for the conic-type submonoid $C=\Z_{\geq 0}\al\op \Z_{\geq 0}\b\leq L$, the associated VOA $V_C$ is $C_1$-cofinite. 
\end{thm}
\begin{proof}	
\begin{figure}[ht]
	\centering
	\begin{tikzpicture}
		\coordinate (Origin)   at (0,0);
		\coordinate (XAxisMin) at (-3,0);
		\coordinate (XAxisMax) at (3,0);
		\coordinate (YAxisMin) at (0,-1);
		\coordinate (YAxisMax) at (0,5);
		\draw [thin, gray,-latex] (XAxisMin) -- (XAxisMax);
		\draw [thin, gray,-latex] (YAxisMin) -- (YAxisMax);
		
		\clip (-8.5,-1) rectangle (3cm,6.4cm); 
		\pgftransformcm{16/25}{0}{-4.16}{sqrt(176/25)}{\pgfpoint{0cm}{0cm}}
		\coordinate (Bone) at (0,2);
		\coordinate (Btwo) at (2,-2);
		\draw[style=help lines,dashed] (-28,-28) grid[step=2cm] (28,28);
		\filldraw[fill=gray, fill opacity=0.3, draw=black] (0,0) rectangle (2,2);
		\foreach \x in {-7,-6,...,14}{
			\foreach \y in {-7,-6,...,7}{
				\node[draw,circle,inner sep=2pt,fill] at (2*\x,2*\y) {};
			}
		}
		\draw [thick,-latex,red] (Origin)
		-- (Bone) node [above] {$\b$};
		\draw [thick,-latex,red] (Origin)
		-- ($(Bone)+(Btwo)$) node [below right] {$\al$};
		\draw [thick,-latex,red] (Origin)
		-- ($2*(Bone)+(Btwo)$) node [above ] {$\b+\al$};
		\draw [thick,-latex,red] (Origin)
		-- ($7*(Bone)+6*(Btwo)$) node [above ] {$\b+(p-1)\al$};
		\draw [thick,-latex,red] (Origin)
		-- ($8*(Bone)+7*(Btwo)$) node [above right] {$\b+p\al$};
	\end{tikzpicture}
	\caption{Fundamental domain of the rank-two lattice \label{fig2}}
\end{figure}

Since $(\cdot|\cdot)$ is positive-definite, we have $(b\b+a\al|b\b+a\al)=2(kb^2-nab+\ell a^2)\geq 0$ for any $a,b\in \Z$, and it is nonzero unless $a=b=0$. In particular, $k(a/b)^2-n(a/b)+\ell\geq 0$ for any rational number $a/b\in \Q$. Since $\Q$ is dense in $\R$, and $n/2k$ is a nonzero rational number, the discriminant $n^2-4k\ell< 0$. Observe that \eqref{eq:conditionforranktwoC_1} is a stronger estimate, hence we need to add it as an assumption. 

Note that $\b$ and $\al$ are in obtuse angle. Among the sequence of vectors $\b,\b+\al,\ds, \b+p\al,\ds$ in $E$, there exists a unique index $p$ such that the angle between $\b+(p-1)\al$ and $\al$ is obtuse, while the angle between $\b+p\al$ and $\al$ is acute or normal. See Figure~\ref{fig2} for an illustration, where the shadow area is the fundamental domain of the lattice $L$. To find such $p$ more explicitly, we note that the interval $[\frac{n}{2\ell}, \frac{n}{2\ell}+1)$ must intersect $\Z$ at a unique positive number. Let 
\begin{equation}\label{eq:defofp}
	p:=\left[\frac{n}{2\ell}, \frac{n}{2\ell}+1\right)\cap \Z\in \Z_{>0}. 
\end{equation}	
Then $\frac{n}{2\ell}+1>p\geq \frac{n}{2\ell}$. In particular, we have
\begin{equation}\label{eq:usefulformula1}
	(\b+(p-1)\al|\al)=-n+2\ell(p-1)<0,\quad (\b+p\al|\al)=-n+2\ell p\geq 0.
\end{equation}
Now consider the following finite subset of the conic-type subVOA $V_C$: $$U=\{\al(-1)\vac, \b(-1)\vac, e^\al, e^{\b}, e^{\b+\al},\ds, e^{\b+p\al}\},$$
where $p$ is given by \eqref{eq:defofp}. Let $W$ be the subspace spanned by elements of the form \eqref{2.13}, with $u^j,u\in U$. From the proof of Lemma~\ref{lm2.14}, we see that to show $V_C=W$, it suffices to show $e^{t\b+s\al}\in W$ for any $t\geq 0$ and $s\geq 0$. 

We use induction on $t\geq 1$ to show $e^{t\b+s\al}\in W$ for all $s\geq 0$. For the base case, note that $(\b+p\al|i\al)=i(-n+2\ell p)\geq 0$ for any $i\geq 0$, in view of \eqref{eq:usefulformula1}. Then $e^{\b+p\al+i\al}\in W$ for all $i\geq 0$, in view of \eqref{2.14}. This shows $e^{\b+s\al}\in W$ for all $s\geq 0$. Assume that for smaller $t$, we have $e^{t\b+s\al}\in W$ for all $s\geq 0$. We want to show $e^{(t+1)\b+s\al}\in W$ for all $s\geq 0$. We divide the range of $s$ into three parts and discuss them one by one. 

Case I. $0\leq s\leq (t+1)(p-1)$. 

For any $0\leq q\leq p-1$, we have $q\leq \frac{n}{2\ell }<\frac{2k}{n}$ as $n< 2\sqrt{k\ell}$. Then $(t\b|\b+q\al)=t(2k-qn)> 0$, and $e^{(t+1)\b+q\al}\in W$ by \eqref{2.14}. For $q=p$, we note that 
\begin{equation}\label{eq:intermediate1}
	(t\b+\al|\b+(p-1)\al)=2kt+(2\ell-tn)(p-1)-n,
\end{equation}
which is non-negative if $2\ell -tn\geq 0$. On the other hand, if $2\ell-tn<0$, by \eqref{eq:usefulformula1} we have 
\begin{align*}
	2kt+(2\ell -tn)(p-1)-n&> 2kt+(2\ell-tn) \left(\frac{n}{2\ell}+1-1\right)-n\\
	&=2kt+n-\frac{tn^2}{2\ell}-n\\
	&=t\cdot \frac{4k\ell-n^2}{2\ell}> 0.
\end{align*}
In other words, the right-hand-side of \eqref{eq:intermediate1} is always non-negative. In fact, we can also see from Figure~\ref{fig2} that the angle between $\b+\al$ and $\b+(p-1)\al$ is acute. 
Since $(t\b+\al)+(\b+(p-1)\al)=(t+1)\b+p\al$, $e^{t\b+\al}\in W$ and $e^{\b+(p-1)\al}\in w$ by the induction hypothesis, we have $e^{(t+1)\b+p\al}\in W$ by \eqref{2.14}.

Furthermore, for any $2\leq j\leq t(p-1)$, it is easy to see from Figure~\ref{fig2} that the angle $\theta$ between $t\b+j\al$ and $\b+(p-1)\al$ is less than the angle between $t\b+\al$ and $\b+(p-1)\al$ which is acute. It follows that $(t\b+j\al|\b+(p-1)\al)\geq 0$ . Since $e^{t\b+j\al}\in W$ and $e^{\b+(p-1)\al}\in W$ by the induction hypothesis, we have $e^{(t+1)\b+(p-1+j)\al}\in W$ for any $2\leq j\leq t(p-1)$ by \eqref{2.14}. This shows $e^{(t+1)\b+ s\al}\in W$ for $0\leq s\leq (t+1)(p-1)$. 

Case II. $(t+1)(p-1)+1 \leq s\leq (t+1)p$.

First, we claim that $(\b+p\al|\b+(p-1)\al)\geq 0$. Indeed, 
\begin{align*}
	(\b+p\al|\b+(p-1)\al)&=2k-(2p-1)n+p(p-1)2\ell\\
	&=2\ell p^2-2(\ell+n)p+(2k+n).
\end{align*}
Regard the right-hand-side as a quadratic form in $p$, then its discriminant $$\Delta=4(n+\ell)^2-8\ell(2k+n)=4(n^2+\ell^2-4\ell k)\leq 0,$$
in view of the condition \eqref{eq:conditionforranktwoC_1}. Hence the quadratic form in $p$ is always non-negative. Note that if we do not require $\Delta\leq 0$, in other words, if $n^2+\ell^2-4\ell k>0$, then it could happen that $p=p_0=\frac{\ell+n}{2\ell}\in \Z$, which satisfies the condition \eqref{eq:defofp}, but $	(\b+p_0\al|\b+(p_0-1)\al)<0$. 

Now, let $a$ be an integer such that $1\leq a\leq t$. By our claim above, 
$$(a\b+a(p-1)\al| (t+1-a)\b+(t+1-a)p\al)=a(t+1-a)\cdot (\b+(p-1)\al|\b+p\al)\geq 0.$$
Then by the induction hypothesis and \eqref{2.14}, we have 
\begin{align*}
	&\exp(a\b+a(p-1)\al+ (t+1-a)\b+(t+1-a)p\al)\\
	&=\exp((t+1)\b+\left(a(p-1)+(t+1-a)p\right)\al)\\
	&=\exp((t+1)\b+((t+1)p-a)\al)\in W. 
\end{align*}
Thus, $e^{(t+1)\b+s\al}\in W$ for any $(t+1)p-t\leq s\leq (t+1)p-1$. Note that $(t+1)p-t=(t+1)(p-1)+1$ which is the lower bound for $s$ in this case.	 Furthermore, since $(t\b+tp\al|\b+p\al)=t\|\b+p\al\|^2> 0$, we have $e^{(t+1)\b+(t+1)p\al}\in W$. 

Case III. $s\geq (t+1)p+1$.

In the previous case, we proved $e^{(t+1)\b+(t+1)p\al}\in W$. For any $i\geq 1$, we have $((t+1)\b+(t+1)p\al| i\al)=2(t+1)pi\ell>0$. Hence $e^{(t+1)\b+((t+1)p+i)\al}\in W$ by \eqref{2.14}. This proves Case III, which implies that $e^{(t+1)\b+s\al}\in W$ for all $s\geq 0$. Now the induction step is complete.   
\end{proof}

\begin{remark}
In Theorem~\ref{thm:ranktwocoincC1}, if $\ell=1$ then the condition~\ref{eq:conditionforranktwoC_1} is satisfied for any $k\geq 1$. Indeed, we know from the proof that $n^2-4k\ell=n^2-4k<0$. Hence $n^2+1-4k\leq 0$. 

In particular, for any rank-two positive-definite even lattice $L$ with a $\Z$-basis $\{\al,\b\}$, if the short basis element has square length $2$ then the conic-type subVOA $V_C=V_{\Z_{\geq 0}\al\op \Z_{\geq 0}\b}$ is $C_1$-cofinite. 
\end{remark}

\subsection{(non) $C_1$-cofiniteness of parabolic-type subVOAs}\label{Sec:5.2}
Unlike the conic-type subVOAs, the parabolic-type subVOAs need not be $C_1$-cofinite. We focus on two concrete examples in this subsection. One of them is $C_1$-cofinite, while the other one is not. 

\begin{prop}\label{prop:C1ofVP}
Let $P=\Z\al\op \Z_{\geq 0}\b$ be the parabolic-type submonoid of the root lattice $A_2$ as in Example~\ref{ex:ranktwoparabolic}. Then $V_P$ is $C_1$-cofinite.  
\end{prop}
\begin{proof}
The proof is similar to the proof of Theorem~\ref{thm:ranktwocoincC1} but is much simpler, we briefly sketch it. Consider the following finite subset:
$$U=\{\al(-1)\vac,\b(-1)\vac, e^\al, e^{-\al}, e^\b, e^{\al+\b}\}.$$
Let $W$ be the subspace spanned by elements of the form \eqref{2.13}, with $u^j,u\in U$. We want to show $W=V_P$. Again, we only need to show $e^{t\b+s\al}\in W$ for any $t\geq 0$ and $s\in \Z$. 

We use induction on $t$. When $t=0$, we have $(\pm \al|\pm i\al)\geq 0$, and so $e^{\pm (i+1)\al}\in W$ for any $i\geq 0$. Hence $e^{s\al}\in  W$ for any $s\in \Z$. Assume the conclusion holds for smaller $t$. Note that for any $m\geq 0$, 
$$(t\b+m\al| \b+\al)=2t-(t+m)+2m=t+m\geq 0.$$
Then by the induction hypothesis and \eqref{2.14}, we have $e^{(t+1)\b+(m+1)\al}\in W$ for any $m\geq 0$. On the other hand, for any $i\geq 0$, we have $(t\b|\b-i\al)=2t+ti>0$, then by the induction hypothesis and \eqref{2.14} again, we have $e^{(t+1)\b-i\al}\in W$ for any $i\geq 0$. Thus, $e^{(t+1)\b+s\al}\in W$ for any $s\in \Z$. The induction step is complete. 
\end{proof}

\begin{remark}
However, the parabolic-type subVOA $V_P=V_{\Z\al\op \Z_{\geq 0}\b}$ is {\bf not} $C_2$-cofinite. The non-$C_2$-cofiniteness of this subVOA is also indicated by the irrationality in Proposition~\ref{prop2.7}. We conjecture that all parabolic-type subVOAs of the form $V_{\Z\al_1\op\ds \op\Z\al_{r-1}\op \Z_{\geq 0}\al_{r}}$ are $C_1$-cofinite.
\end{remark}

Recall that Borel-type subVOAs are of parabolic-type by definition. Next, we want to show that the rank-two Borel-type subVOA $V_B$ in Example~\ref{ex2.7} is not $C_1$-cofinite. Thus, there are no uniform results on the $C_1$-cofiniteness for general parabolic-type subVOAs.

\begin{lm}\label{lm:nonC_1forVB}
Let $B=\{ m\al+n\b: m,n\in \Z, n>0\}+ \Z_{\geq 0}\al$ be the Borel-type submonoid of $A_2$ as in Example~\ref{ex:ranktwoparabolic}. Then $e^{\b-m\al}\notin C_1(V_B)$ for any $m\geq 1$.
\end{lm}
\begin{proof}
By the definition of $C_1(V)$ \eqref{eq:defofC1} and $V_B$, we have 
\begin{align*}
	C_1(V_B)&=\spn\{u_{-1}v: u\in M_{\hat{\h}}(1,\ga)\cap (V_B)_+, v\in M_{\hat{\h}}(1,\ga')\cap(V_B)_+, \ga,\ga'\in B \}\\
	&\ +\spn\{L(-1)w:w\in M_{\hat{\h}}(1,\theta)\cap(V_B)_+, \theta\in B \}.
\end{align*}
Suppose there exists some $m\geq 1$ such that $e^{\b-m\al}\in C_1(V_B)$. Since the Heisenberg modules $M_{\hat{\h}}(1,\ga)$ are in direct sum in $V_B$, and they are invariant under $L(-1)$, then by \eqref{2.8} we have $\b-m\al\in B+B$. Write $\b-m\al=\ga+\ga'$ for some $\ga,\ga'\in B$. By the definition of $B$, the only possibilities for the pair $(\ga,\ga')$ are 
$$(\ga,\ga')=(\b-(m+k)\al,k\al)\quad \mathrm{or}\quad (\ga,\ga')=(k\al,\b-(m+k)\al),$$
where $k\geq 0$. Consider the following subspaces of $C_1(V_B)$:
\begin{align*}
	W_1:&=\spn\{u_{-1}v: u\in M_{\hat{\h}}(1,\b-(m+k)\al)\cap (V_B)_+,\ v\in M_{\hat{\h}}(1,k\al)\cap (V_B)_+ \},\\
	W_2:&=\spn\{u'_{-1}v': u'\in M_{\hat{\h}}(1,k\al)\cap (V_B)_+,\ v'\in M_{\hat{\h}}(1,\b-(m+k)\al)\cap (V_B)_+ \},\\
	W_3:&=\spn\{L(-1)w: w\in  M_{\hat{\h}}(1,\b-m\al)\cap (V_B)_+\}.
\end{align*}
We have $e^{\b-m\al}\in W_1+W_2+W_3$.

Note that $\wt (e^{\b-m\al})=(\b-m\al|\b-m\al)/2=m^2+m+1$. For any homogeneous elements $u\in M_{\hat{\h}}(1,\b-(m+k)\al)\cap (V_B)_+ $ and $ v\in M_{\hat{\h}}(1,k\al)\cap (V_B)_+$, if $k>0$ then 
\begin{align*}
	\wt (u_{-1}v)&=\wt u+\wt v\geq \frac{(\b-(m+k)\al|\b-(m+k)\al)}{2}+\frac{(k\al|k\al)}{2}\\
	&=1+(m+k)+(m+k)^2+k^2> \wt (e^{\b-m\al}).
\end{align*}
If $k=0$, then $v\in M_{\hat{\h}}(1,0)\cap (V_B)_+\ssq \sum_{t=1}^\infty M_{\hat{\h}}(1,0)_t$, and 
$$\wt (u_{-1}v)\geq \frac{(\b-m\al|\b-m\al)}{2}+\wt (v)\geq (m^2-m+1)+1> \wt(e^{\b-m\al}).$$
Thus, $\wt (u_{-1}v)>\wt (e^{\b-m\al})$. Similarly, we can show $\wt (u'_{-1}v')>\wt (e^{\b-m\al})$ for any $u'_{-1}v'\in W_2$. Finally, for any $L(-1)w\in W_3$, with $w\in  M_{\hat{\h}}(1,\b-m\al)\cap (V_B)_+$, we have $\wt (L(-1)w)=\wt(w)+1> \wt(e^{\b-m\al})$ since $\wt (w)\geq \wt(e^{\b-m\al})$. Therefore, for any homogeneous element $x\in W_1+W_2+W_3$, we have $\wt (x)>\wt(e^{\b-m\al})$. This contradicts $e^{\b-m\al}\in W_1+W_2+W_3$. 
\end{proof}

\begin{prop}\label{prop:VBnonC1}
Let $B=\{ m\al+n\b: m,n\in \Z, n>0\}+ \Z_{\geq 0}\al$ be the Borel-type submonoid of $A_2$ as in Example~\ref{ex:ranktwoparabolic}. Then $V_B$ is {\bf not} $C_1$-cofinite.  
\end{prop}
\begin{proof}
It suffices to show elements in the set $S=\{\vac+C_1(V_B), e^{\b-m\al}+C_1(V_B):m\geq 1 \}$ are linearly independent in $V_B/C_1(V_B)$. Indeed, by Lemma~\ref{lm:nonC_1forVB}, these elements are nonzero. Since $L(0)C_1(V)\ssq C_1(V)$, it induces a well-defined linear map $L(0): V_B/C_1(V_B)\ra V_B/C_1(V_B)$. For  $m=1,2,\dots$, we have
$$L(0)(e^{\b-m\al}+C_1(V_B))=(m^2+m+1)(e^{\b-m\al}+C_1(V_B)).$$
i.e., the elements in $S$ are in different grading subspaces of $V_B/C_1(V_B)$. Hence the elements in $S$ are linearly independent. 
\end{proof}

\begin{remark}
One can take a submonoid $M=\{m\al+n\b: m\geq n\geq 1\}\cup \{0\}$ of $B$ as in Example~\ref{ex:ranktwoparabolic} and construct another example $V_M$ of CFT-type non-$C_1$-cofinite VOA. This particular VOA shows that the Noetherianity of Zhu's algebra $A(V)$ does not hold if $V$ is not $C_1$-cofinite. See \cite[Section 4]{Liu25} for more details.
\end{remark}

\section{The rank-one Borel-type subVOA $V_B$ of $V_{\Z\al}$}\label{sec3}	

In this Section, we fix the rank-one lattice $L=\Z\al$, with $(\al|\al)=2N$ for some $N\geq 1$, and $\epsilon(m\al,n\al)=1$, for all $m,n\in \Z$, and study its Borel-type subVOA $V_{B}$, where $B=\Z_{\geq 0}\al$. Note that $V_B$ is also of conic-type and parabolic-type, see Example~\ref{ex:rankoneBorel}. 
Our goal is to show that the Zhu's algebra $A(V_B)$ has the following identification as associative algebras: 
\begin{equation}\label{eq:6.1}
A(V_B)\cong \C\<x,y\>/\<y^2,yx+Ny,xy-Ny\>,
\end{equation}
where $\C\<x,y\>$ is the tensor algebra on generators $x$ and $y$.

First, we recall the definition and basic properties of Zhu's algebra $A(V)$ in \cite{Z,FZ}. Let $(V,Y,\vac,\om)$ be a VOA. For homogeneous elements $a,b\in V$, define
\begin{align}
a\circ b&:=\Res_{z} Y(a,z)b\frac{(1+z)^{\wt a}}{z^2}=\sum_{j\geq 0} \binom{\wt a}{j} a_{j-2}b\label{3.1},\\
a\ast b&:=\Res_{z}Y(a,z)b\frac{(1+z)^{\wt a}}{z}=\sum_{j\geq 0} \binom{\wt a}{j} a_{j-1}b\label{3.2}.
\end{align}
Let $O(V)=\spn\{a\circ b:a,b\in V\}$, and let $A(V):=V/O(V)$. By \cite[Theorem 2.1.1]{Z}, $O(V)$ is a two-sided ideal with respect to $\ast$, and $A(V)$ is an associative algebra with respect to $\ast$, with the unit element $[\vac]$. By \cite[Lemma 2.1.3]{Z}, we also have the following formulas:
\begin{align}
a\ast b&\equiv \Res_{z} Y(b,z)a\frac{(1+z)^{\wt b-1}}{z} \pmod{O(V)},\label{3.3}\\
a\ast b-b\ast a&\equiv \Res_{z}Y(a,z)b(1+z)^{\wt a-1}\pmod{O(V)},\label{3.4}
\end{align}
for any homogeneous $a,b\in V$. Furthermore, if $m\geq n\geq 0$, one has 
\begin{equation}\label{3.5}
\Res_{z} Y(a,z)b\frac{(1+z)^{\wt a+n}}{z^{2+m}}\equiv 0\pmod{O(V)}.
\end{equation}

\subsection{A spanning set of $O(V_B)$}\label{sec:6.1} To determine the structure of $A(V_B)$, we first give a description of $O(V_B)$ in its definition. 

\begin{df}
Let $V_B=V_{\Z_{\geq 0}\al}$ be the Borel-type subVOA of $V_{\Z\al}$, with $(\al|\al)=2N$.
Let $O'$ be the subspace of $V_B$ spanned by the following elements:
\begin{equation}\label{3.8}
	\begin{cases}
		&\al(-n-2)u+\al(-n-1)u, \qquad u\in V_B,\ \mathrm{and}\ n\geq 0,\\
		&\al(-1)v+Nv,\qquad v\in \bigoplus_{m\geq 1} M_{\hat{\h}}(1,m\al),\\
		&M_{\hat{\h}}(1,k\al),\qquad k\geq 2.
	\end{cases}
\end{equation}
\end{df}

We want to show that $O(V_B)=O'$. Then the structure of $A(V_B)$ would be clear. 

\subsubsection{Proof of $O'\ssq O(V_B)$}
It is not difficult to show $O'\ssq O(V_B)$. Indeed, it is clear from \eqref{3.5} that $\al(-n-2)u+\al(-n-1)u\in O(V_B)$, for all $u\in V_B$ and $n\geq 0$. By \cite[Theorem 2.1.1]{Z}, we also have 
\begin{equation}\label{3.9}
a\ast O(V_B) \subset O(V_B),\quad \mathrm{and}\quad O(V_B)\ast a\subset O(V_B),\quad a\in V.
\end{equation}

\begin{lm}\label{lm3.2}
For any $k\geq 2$, we have $M_{\hat{\h}}(1,k\al)\subset O(V_B)$.
\end{lm}
\begin{proof}
By the proof of Proposition \ref{prop3.1}, we have $e^{k\al}\in O(V_B)$, for any $k\geq 2$. By \eqref{3.3}, we have $u\ast \al(-1)\vac\equiv \al(-1)u\pmod{O(V_B)}.$ Now by \eqref{3.5} and \eqref{3.9}, we have
\begin{align*}
	\al(-n_{1}-1)\dots \al(-n_r-1)&e^{k\al}\equiv (-1)^{n_{1}+\dots +n_r} \al(-1)^re^{k\al}\pmod{O(V_B)}\\
	&\equiv (-1)^{n_{1}+\dots +n_r} e^{k\al} \ast (\al(-1)\vac)\ast\dots \ast(\al(-1)\vac)\pmod{O(V_B)}\\
	&\equiv 0\pmod{O(V_B)},
\end{align*}
for any $k\geq 2$, $r\geq 1$, and $n_{1},\dots n_{r}\geq 0$, where the last congruence follows from \eqref{3.9}. Thus we have $M_{\hat{\h}}(1,k\al)\subset O(V_B)$, for all $k\geq 2$. 
\end{proof}

\begin{lm}\label{lm3.3}
For any $v\in \bigoplus_{m\geq 1} M_{\hat{\h}}(1,m\al)$, we have $\al(-1)v+Nv\in O(V_B)$. 
\end{lm}
\begin{proof}
If $m\geq 2$ and $v\in M_{\hat{\h}}(1,m\al)$, then by Lemma \ref{lm3.2}, we have $v\in O(V_B)$, and $$\al(-1)v+Nv\equiv v\ast (\al(-1)\vac)+Nv\equiv 0\pmod{O(V_B)},$$
by \eqref{3.9}. Now let $m=1$, by the proof of Proposition \ref{prop3.1}, we have $\al(-1)e^\al+Ne^\al=e^\al\circ \vac\equiv 0\pmod{O(V_B)}.$ Let $v=\al(-n_{1}-1)\dots \al(-n_r-1)e^{\al}$ be a general spanning element of $M_{\hat{\h}}(1,\al)$, where $r\geq 1$, and $n_{1},\dots n_{r}\geq 0$. Since $[\al(-1),\al(-p)]=0$ for all $p\geq 1$, we have
\begin{align*}
	\al(-1) v+Nv&= \al(-n_{1}-1)\dots \al(-n_r-1)(\al(-1)e^\al+Ne^\al)\\
	&\equiv  (-1)^{n_{1}+\dots +n_r} \al(-1)^r(\al(-1)e^\al+Ne^\al)\pmod{O(V_B)}\\
	&\equiv  (-1)^{n_{1}+\dots +n_r} (\al(-1)e^{\al}+Ne^\al) \ast (\al(-1)\vac)\ast\dots \ast(\al(-1)\vac)\pmod{O(V_B)}\\
	&\equiv 0\pmod{O(V_B)},
\end{align*}
where the last congruence follows from \eqref{3.9} and the fact that $\al(-1)e^{\al}+Ne^\al\in O(V_B)$. 
\end{proof}

It follows from Lemma \ref{3.2} and Lemma \ref{lm3.3} that $O'\ssq O(V_B)$.

\subsubsection{Outline of the strategy for proving $O(V_B)\ssq O'$}


Conversely, by the definition of $O(V_B)$ \eqref{3.1}, to show $O(V_B)\ssq O'$, we need to show that $a\circ u=\Res_{z}Y(a,z)u((1+z)^{\wt a}/z^2)\in O'$, for any homogeneous $a,u\in V_B$.

Note that if $a\in M_{\hat{\h}}(1,m\al)$ and $u\in M_{\hat{\h}}(1,n\al)$ for some $m,n\geq 1$, then by \eqref{2.3} and \eqref{2.8},
$$a\circ u=\Res_{z}Y(a,z)u\frac{(1+z)^{\wt a}}{z^2}\in M_{\hat{\h}}(1,(m+n)\al)((z))\subset O'((z)),$$
as $m+n\geq 2$ and $M_{\hat{\h}}(1,k\al)\subset O'$ for any $k\geq 2$. Since $V_B=M_{\hat{\h}}(1,0)\op \bigoplus_{m=1}^\infty M_{\hat{\h}}(1,m\al)$ as a vector space, it suffices to show:
\begin{equation}\label{3.10}
a\circ u\in O',\quad \mathrm{for}\quad \begin{cases}a\in M_{\hat{\h}}(1,\al) &\mathrm{and}\ \ \ u\in M_{\hat{\h}}(1,0),\\ \mathrm{or}&\\ a\in M_{\hat{\h}}(1,0) &\mathrm{and}\ \ \ u\in M_{\hat{\h}}(1,\al).\end{cases}
\end{equation}

We will only present a detailed proof for the first case of \eqref{3.10}. i.e., $M_{\hat{\h}}(1,\al)\circ M_{\hat{\h}}(1,0)\ssq O'$. The second case follows from a similar argument (see also (3.1.5) and (3.1.6) in \cite{FZ}), we omit the detail. Our strategy of proving  $M_{\hat{\h}}(1,\al)\circ M_{\hat{\h}}(1,0)\ssq O'$ can be outlined as follows: 
\begin{enumerate}[label=Step \arabic*.,align=left]
\item \label{step1} Let $u=a(-m_1)\ds a(-m_s)\vac\in M_{\hat{\h}}(1,0)$. Use an induction on the length $s$ of $u$ to show that 
$$\Res_{z} Y(e^\al,z)u\frac{(1+z)^N}{z^{2+n}}\in O',\quad \mathrm{for\ any}\ n\geq 0.$$
Note that $N=\wt (e^\al)$. 
\item \label{step2} Let $a=\al(-k)e^\al$ for some $k\geq 1$. Use an induction on the number $k\geq 1$, together with the result in Step 1, to show that for any $u\in M_{\hat{\h}}(1,0)$, 
$$\Res_z Y(a(-k)e^\al,z)u \frac{(1+z)^{N+k}}{z^{2+n}} \in O',\quad \mathrm{for\ any}\ n\geq 0.$$
Note that $N+k=\wt (a(-k)e^\al)$. 
\item \label{step3} Let $a=\al(-k)\al(-n_2)\dots \al(-n_r)e^\al\in M_{\hat{\h}}(1,\al)$, where $k\geq 1$. Use an induction on the length $r$ of $a$ (with the base case given by Steps 1 and 2), together with another induction on the number $k\geq 1$, to show that for any $u\in M_{\hat{\h}}(1,0)$, 
$$	\Res_{z} Y(\al(-k)\al(-n_{2})\ds \al(-n_{r})e^\al,z)u \frac{(1+z)^{N+k+n_{2}+\ds+n_{r}}}{z^{2+n}}\in O',\quad \mathrm{for\ any}\ n\geq 0.$$
Note that $N+k+n_{2}+\ds+n_{r}=\wt (\al(-k)\al(-n_{2})\ds \al(-n_{r})e^\al)$. 
\end{enumerate} 
Then the conclusion in Step 3 would indicate $M_{\hat{\h}}(1,\al)\circ M_{\hat{\h}}(1,0)\ssq O'$, and the proof of the  inclusion $O(V_B)\ssq O'$ would be complete. 
\subsection{Proof of the main Theorem}\label{sec3.2}

Now we carry out the details for proving: 
\begin{equation}\label{eq:main}
a\circ u\in O',\quad \mathrm{for}\quad \begin{cases}
	a=\al(-n_1)\dots \al(-n_r)e^\al\in M_{\hat{\h}}(1,\al),\\
	u=\al(-m_1)\ds \al(-m_s)\vac\in M_{\hat{\h}}(1,0),
\end{cases}
\end{equation}
where $r,s\geq 0$, $n_1\geq\ds \geq n_r\geq 1$, and $m_1\geq \ds \geq m_s\geq 1$, following  the steps outlined above. 
Then by our discussions in Section~\ref{sec:6.1}, this would lead to a proof of the following main theorem in this Section: 

\begin{thm}\label{thm3.8}
Let $V_B=V_{\Z_{\geq 0}\al}$ be the Borel-type subVOA of $V_{\Z\al}$, with $(\al|\al)=2N$. Then 
\begin{equation}\label{3.22}
	\begin{split}O(V_B)=O'=\spn \bigg\{ \al(-n-2)u+\al(-n-1)u,\ \al(-1)v+v,\ M_{\hat{\h}}(1,k\al):&
		\\n\geq 0,\  u\in V_B, \ v\in \bigoplus_{m\geq 1}M_{\hat{\h}}(1,m\al),\ k\geq 2\bigg\}.
	\end{split}
\end{equation}
\end{thm}

\subsubsection{Proof of the claim in \ref{step1}}
We prove a Lemma first. 
\begin{lm}\label{lm3.4}
For any $m\geq 1$, we have $\al(-m)O'\subset O'$. For any $a\in M_{\hat{\h}}(1,\al)$, we have $L(-1)a+L(0)a\in O'$. 
\end{lm}
\begin{proof}
Since $[\al(-m),\al(-n)]=0$ for any $m,n\geq 1$, and $\al(-m) M_{\hat{\h}}(1,k\al)\subset M_{\hat{\h}}(1,k\al)$, for any $k\geq 0$, we have $\al(-m)O'\subset O'$, in view of \eqref{3.8}. 

Let $a=\al(-m_1)\dots \al(-m_s)e^\al\in M_{\hat{\h}}(1,\al)$, where $r\geq 0$ and $m_1,\ds,m_s\geq 1$. Since $L(-1)e^\al=(e^\al)_{-2}\vac=\al(-1)e^\al$, and $[L(-1),\al(-m)]=ma(-m-1)$, we have
\begin{align*}
	&L(-1)\al(-m_1)\dots \al(-m_r)e^\al+L(0)\al(-m_1)\dots \al(-m_s)e^\al\\
	&=\al(-m_1)\dots \al(-m_r)\al(-1)e^\al+\sum_{j=1}^s m_j\cdot\al(-m_1)\dots \al(-m_j-1)\dots \al(-m_s)e^\al\\
	&\quad+(m_{1}+\dots +m_k+N) \al(-m_1)\dots \al(-m_s)e^\al\\
	&=\al(-m_1)\dots \al(-m_s)(\al(-1)e^\al+Ne^\al)\\
	&\quad+\sum_{j=1}^s (\al(-m_j-1)+\al(-m_j)) \al(-m_1)\dots \widehat{\al(-m_j)} \dots \al(-m_s)e^\al\\
	&\equiv 0\pmod{O'},
\end{align*}
where the last congruence follows from $\al(-1)e^\al+Ne^\al\in O'$, $\al(-m)O'\subset O'$ for any $m\geq 1$, and $\al(-n-1)v+\al(-n)v\in O'$ for all $v\in M_{\hat{\h}}(1,\al)$ and $n\geq 1$, in view of \eqref{3.8}. 
\end{proof}

\begin{prop}\label{prop3.5}
Let $u\in M_{\hat{\h}}(1,0)$, and $n\geq 0$. We have 
\begin{equation}\label{3.11}
	\Res_{z} Y(e^\al,z)u\frac{(1+z)^N}{z^{2+n}}\in O'.
\end{equation}
\end{prop}
\begin{proof}
We use induction on the length $s$ of a spanning element $u= \al(-m_{1})\dots \al(-m_s)\vac$ of $M_{\hat{\h}}(1,0)$, where $m_{1},\dots ,m_{s}\geq 1$. The base case is $u=\vac$. First, we note that $$e^\al_{-j-1}\vac=\frac{1}{j!} (L(-1)^je^\al)_{-1}\vac=\frac{1}{j!} L(-1)^j e^\al,\quad j\geq 0.$$  Since $e^\al_{-j-1}\vac\in M_{\hat{\h}}(1,\al)$ for any $j\geq 0$, by Lemma \ref{lm3.4} we have 
\begin{align*}
	L(-1)^j e^\al&\equiv -L(0)L(-1)^{j-1}e^\al\pmod{O'}\\
	&=(-1) (N+j-1) L(-1)^{j-1} e^\al \\
	&\vdots\\
	&\equiv (-1)^j (N+j-1)(N+j-2)\ds (N+1)N e^\al \pmod{O'}.
\end{align*}
Then it follows from the definition of binomial coefficients that
\begin{align*}
	Y(e^\al,z)\vac&=\sum_{j\geq 0}e^\al_{-j-1}\vac z^j=\sum_{j\geq 0} \frac{1}{j!} L(-1)^j  z^je^\al\\
	&\equiv \sum_{j\geq 0} (-1)^j\frac{(N+j-1)(N+j-2)\ds (N+1)N}{j!} z^je^\al\pmod{O'}\\
	&=\sum_{j\geq 0} \frac{(-N-j+1)(-N-j+2)\ds (-N-1)(-N)}{j!}z^j e^\al\numberthis \label{3.12}\\
	&=\sum_{j\geq 0} \binom{-N}{j}z^je^\al=(1+z)^{-N} e^\al .
\end{align*}
Now by \eqref{3.11} and \eqref{3.12}, and the assumption that $n\geq 0$, we have
$$
\Res_{z} Y(e^\al,z)\vac\frac{(1+z)^N}{z^{2+n}}\equiv \Res_{z} (1+z)^{-N} \frac{(1+z)^N}{z^{2+n}}e^\al=\Res_{z}\frac{1}{z^{2+n}}e^\al=0\pmod{O'}.
$$
This finishes the proof of the base case. Assume the conclusion holds for smaller length $s$. Note that for any $m\geq 1$, we have
$$
[\al(-m),Y(e^\al,z)]=\sum_{i\geq 0} \binom{-m}{i}Y(\al(i)e^\al,z)z^{-m-i}=2NY(e^\al,z) z^{-m}.
$$
Then by $\al(-m)O'\subset O'$ in Lemma \ref{lm3.4}, the base case and the induction hypothesis, we have
\begin{align*}
	&\Res_{z} Y(e^\al,z) \al(-m_1)\ds \al(-m_s)\vac\frac{(1+z)^N}{z^{2+n}}\\
	&=\Res_{z} \al(-m_1)\ds \al(-m_s)Y(e^\al,z)\vac\frac{(1+z)^N}{z^{2+n}}\\
	&\quad-\sum_{j=1}^s \Res_{z} \al(-m_1) \ds \al(-m_{j-1})[\al(-m_j),Y(e^\al,z)]\al(m_{j+1})\ds \al(-m_s) \vac\frac{(1+z)^N}{z^{2+n}}\\
	&\equiv -\sum_{j=1}^s 2N\Res_{z} \al(-m_1)\ds \al(-m_{j-1})Y(e^\al,z)\al(-m_{j+1})\ds \al(-m_s)\vac \frac{(1+z)^N}{z^{2+n+m_j}}\pmod{O'}\\
	&\equiv 0\pmod{O'}
\end{align*}
where the last congruence follows from the induction hypothesis, which indicates that $$\Res_{z} Y(e^\al,z)\al(-m_{j+1})\ds \al(-m_s)\vac ((1+z)^N)/z^{2+n+m_j}\in O'.$$ Hence \eqref{3.11} holds for any $u\in M_{\hat{\h}}(1,0)$ and $n\geq 0$.
\end{proof}

\subsubsection{Proof of the claim in \ref{step2}}

\begin{prop}\label{lm3.6}
For any $k\geq 1$, $n\geq 0$, and $u\in M_{\hat{\h}}(1,0)$, we have
\begin{equation}\label{3.14}
	\Res_{z} Y(\al(-k)e^\al,z)u\frac{(1+z)^{N+k}}{z^{2+n}}\in O'.
\end{equation}
\end{prop}
\begin{proof}
It is easy to derive the following formula from the Jacobi identity of VOAs: 
\begin{equation}\label{3.15}
	Y(\al(-1)v,z)=\sum_{j\geq 0} \al(-j-1) Y(v,z)z^{j}+\sum_{j\geq 0} Y(v,z) \al(j) z^{-j-1},\quad v\in V_{\Z\al}.
\end{equation}
Now we prove \eqref{3.14} by induction on $k$. When $k=1$, by \eqref{3.15} we have
\begin{align*}
	&\Res_{z} Y(\al(-1)e^\al,z)u\frac{(1+z)^{N+1}}{z^{2+n}}\\
	&=\Res_{z} \sum_{j\geq 0} \al(-j-1)Y(e^\al,z)u z^j\frac{(1+z)^{N+1}}{z^{2+n}}+\Res_{z}\sum_{j\geq 0} Y(e^\al,z)\al(j)u\frac{(1+z)^{N+1}}{z^{2+n+j+1}}\\
	&= \Res_{z} \left(\sum_{j\geq 0} \al(-j-1)Y(e^\al,z)u z^j\frac{(1+z)^N}{z^{2+n}} +\sum_{j\geq 0} \al(-j-1)Y(e^\al,z)u z^{j+1}\frac{(1+z)^N}{z^{2+n}}\right)\\
	&\quad+\Res_{z}\sum_{j\geq 0} Y(e^\al,z)\al(j)u\frac{(1+z)^{N}}{z^{2+n+j+1}}+ \Res_{z}\sum_{j\geq 0} Y(e^\al,z)\al(j)u\frac{(1+z)^{N}}{z^{2+n+j}}\\
	&\equiv \Res_{z} \al(-1)Y(e^\al,z)u\frac{(1+z)^N}{z^{2+n}} \\
	&\quad+\sum_{j\geq 0}(\al(-j-2)+\al(-j-1))\Res_{z} Y(e^\al,z) uz^{j+1}\frac{(1+z)^N}{z^{2+n}}+0\pmod{O'}\\
	&\equiv \Res_{z} \al(-1)Y(e^\al,z)u\frac{(1+z)^N}{z^{2+n}}\pmod{O'},
\end{align*}
where the first congruence follows from Proposition \ref{prop3.5}, as $n+j\geq 0$, and the second congruence follows from \eqref{3.8}. Furthermore, by Proposition \ref{prop3.5} again, we have
\begin{align*}
	&\Res_{z}\al(-1)Y(e^\al,z)u\frac{(1+z)^N}{z^{2+n}}\\
	&=\Res_{z}\left( Y(e^\al,z)\al(-1)u\frac{(1+z)^N}{z^{2+n}}+\sum_{j\geq 0} \binom{-1}{j} z^{-1-j} Y(\al(j)e^\al,z)u \frac{(1+z)^N}{z^{2+n}}\right)\\
	&\equiv 0+ \Res_{z} 2NY(e^\al,z)u\frac{(1+z)^N}{z^{2+n+1}}\pmod{O'}\\
	&\equiv 0\pmod{O'}.
\end{align*}
This proves \eqref{3.14} for $k=1$. Assume \eqref{3.14} holds for smaller $k$. Note that  $[L(-1),\al(-k)]=\frac{1}{k}\al(-k-1)$ and $\al(-k)L(-1)e^\al=\al(-k)\sum_{i\geq 0} \frac{\al}{\sqrt{2N}}(-1-i)\frac{\al}{\sqrt{2N}}(i)e^\al=\al(-1)\al(-k)e^\al$. By \eqref{3.15}, 
\begin{align*}
	&\Res_{z} Y(\al(-k-1)e^\al,z)u\frac{(1+z)^{N+k+1}}{z^{2+n}}\\
	&= \frac{1}{k}\Res_{z}\left(Y(L(-1)\al(-k)e^\al,z)u\frac{(1+z)^{N+k+1}}{z^{2+n}}-  Y(\al(-k)L(-1)e^\al,z)u\frac{(1+z)^{N+k+1}}{z^{2+n}}\right)\\
	&=- \frac{1}{k}\Res_{z}  \left(Y(\al(-k)e^\al,z)u\frac{d}{dz}\left(\frac{(1+z)^{N+k+1}}{z^{2+n}}\right)+ Y(\al(-1)\al(-k)e^\al,z)u\frac{(1+z)^{N+k+1}}{z^{2+n}}\right)\\
	&=-\frac{N+k+1}{k} \Res_{z} Y(\al(-k)e^\al,z)u\frac{(1+z)^{N+k}}{z^{2+n}}\\
	&\quad+\frac{2+n}{k}\Res_{z} Y(\al(-k)e^\al,z)u\frac{(1+z)^{N+k}}{z^{3+n}}+\frac{2+n}{k}\Res_{z} Y(\al(-k)e^\al,z)u\frac{(1+z)^{N+k}}{z^{2+n}}\\
	&\quad-\frac{1}{k}\Res_{z} \sum_{j\geq 0} \al(-j-1) Y(\al(-k)e^\al,z)u z^j\frac{(1+z)^{N+k+1}}{z^{2+n}}\\
	&\quad-\frac{1}{k} \Res_{z}\sum_{j\geq 0}Y(\al(-k)e^\al,z)\al(j)u \frac{(1+z)^{N+k+1}}{z^{2+n+j+1}}\\
	&\equiv 0-\frac{1}{k}\Res_{z}\sum_{j\geq 0} \al(-j-1)\left(Y(\al(-k)e^\al,z)uz^j+Y(\al(-k)e^\al,z)uz^{j+1}\right)\frac{(1+z)^{N+k}}{z^{2+n}} \\
	&\quad-\frac{1}{k} \Res_{z} \sum_{j\geq 0}Y(\al(-k)e^\al,z)\al(j)u \frac{(1+z)^{N+k}(1+z)}{z^{2+n+1+j}}\pmod{O'}\\
	&=	-\frac{1}{k} \Res_{z} \al(-1)Y(\al(-k)e^\al,z)u\frac{(1+z)^{N+k}}{z^{2+n}}\\
	&\quad-\frac{1}{k} \Res_{z} \sum_{j\geq 0}\left(\al(-j-2)+\al(-j-1)\right)Y(\al(-k)e^\al,z)u z^{j+1} \frac{(1+z)^{N+k}}{z^{2+n}}\\
	&\quad-\frac{1}{k} \Res_{z} \sum_{j\geq 0}Y(\al(-k)e^\al,z)\al(j)u \frac{(1+z)^{N+k}(1+z)}{z^{2+n+1+j}}\\
	&\equiv -\frac{1}{k} \Res_{z} \al(-1)Y(\al(-k)e^\al,z)u\frac{(1+z)^{N+k}}{z^{2+n}}\pmod{O'},
\end{align*}
where the first congruence follows from the induction hypothesis \eqref{3.14}, and the second congruence follows from \eqref{3.8} and the induction hypothesis. By the Jacobi identity and the Heisenberg relation $[\al(j),\al(-k)]=\delta_{j,k}k K$ for any $j\geq 0$, we have
\begin{align*}
	&-\frac{1}{k} \Res_{z} \al(-1)Y(\al(-k)e^\al,z)u\frac{(1+z)^{N+k}}{z^{2+n}}\\
	&=-\frac{1}{k} \Res_{z}Y(\al(-k)e^\al,z)\al(-1)u\frac{(1+z)^{N+k}}{z^{2+n}}-\frac{1}{k} \Res_{z}\sum_{j\geq 0} \binom{-1}{j} z^{-1-j} Y(\al(j)\al(-k)e^\al,z)u \frac{(1+z)^{N+k}}{z^{2+n}}\\
	&\equiv 0-\frac{1}{k} \Res_{z} (-1)^kkY(e^\al,z)u(1+z)^k\frac{(1+z)^N}{z^{2+n+1+k}}\pmod{O'}\\
	&=-\Res_{z} (-1)^{k}\sum_{i\geq 0} \binom{k}{i} Y(e^\al,z)u\frac{(1+z)^N}{z^{2+n+1+i}}\\
	&\equiv 0\pmod{O'},
\end{align*}
where the first congruence follows from the induction hypothesis, and the second congruence follows from Proposition \ref{prop3.5}. Therefore, we have
$$\Res_{z} Y(\al(-k-1)e^\al,z)u\frac{(1+z)^{N+k+1}}{z^{2+n}}\in O'.$$
So \eqref{3.14} holds for $k+1$, and the inductive step is complete. 
\end{proof}

\subsubsection{Proof of the claim in \ref{step3}}


\begin{prop}\label{prop3.7}
Let $a= \al(-n_{1})\dots \al(-n_{r})e^\al\in M_{\hat{\h}}(1,\al)$, where $r\geq 0$ and $n_{1},\dots,n_{r}\geq 1$. We have 
\begin{equation}\label{3.16}
	\Res_{z} Y(a,z)u\frac{(1+z)^{\wt a}}{z^{2+n}}\in O',
\end{equation}
for any $n\geq 0$ and $u\in M_{\hat{\h}}(1,0)$. In particular, \eqref{eq:main} holds true. 
\end{prop}
\begin{proof}
We use induction on the length $r$ of $a$ to show \eqref{3.16}. By Propositions \ref{prop3.5} and \ref{lm3.6}, \eqref{3.16} holds for $a=e^\al$ or $a=\al(-k)e^\al$ for some $k\geq 1$. Now let $r\geq 2$. The induction hypothesis is the assumption that 
\begin{equation}\label{3.17}
	\Res_{z}Y(\al(-n_{2})\dots \al(-n_{r})e^\al,z)u\frac{(1+z)^{N+n_{2}+\ds +n_{r}}}{z^{2+n}}\in O',
\end{equation}
for $n_{2},\ds,n_{r}\geq 1$, $n\geq 0$, and $u\in M_{\hat{\h}}(1,0)$.	

First, we claim that
\begin{equation}\label{3.18}
	\Res_{z}Y(\al(-1)\al(-n_{2})\ds \al(-n_{r})e^\al,z)u\frac{(1+z)^{N+n_{2}+\ds+n_{r}+1}}{z^{2+n}}\in O'.
\end{equation}
Denote $N+n_{2}+\ds +n_{r}$ by $m$, then $\wt (\al(-n_{2})\ds \al(-n_{r})e^\al)=m$. By \eqref{3.15}, \eqref{3.8}, and the induction hypothesis, we have 
\begin{align*}
	&\Res_{z}Y(a(-1)\al(-n_{2})\ds \al(-n_{r})e^\al,z)u\frac{(1+z)^{1+m}}{z^{2+n}}\\
	&=\Res_{z}\sum_{j\geq 0} a(-j-1)Y(\al(-n_{2})\ds \al(-n_{r})e^\al, z)uz^{j}\frac{(1+z)^{1+m}}{z^{2+n}}\\
	&\quad + \Res_{z}\sum_{j\geq 0} Y(\al(-n_{2})\ds \al(-n_{r})e^\al,z)(a(j)u)\frac{(1+z)^{1+m}}{z^{2+n+j+1}}\\
	&\equiv \Res_{z}\sum_{j\geq 0} a(-j-1)Y(\al(-n_{2})\ds \al(-n_{r})e^\al,z)uz^{j}\frac{(1+z)^{m}}{z^{2+n}}\\
	&\quad +\Res_{z}\sum_{j\geq 0} a(-j-1)Y(a_{2}(\al(-n_{2})\ds \al(-n_{r})e^\al,z)uz^{j+1}\frac{(1+z)^{m}}{z^{2+n}}\pmod{O'}\\
	&=
	\Res_{z}a(-1)Y(\al(-n_{2})\ds \al(-n_{r})e^\al,z)u\frac{(1+z)^{m}}{z^{2+n}}\\
	&\quad +\Res_{z}\sum_{j\geq 0} a(-j-2)Y(\al(-n_{2})\ds \al(-n_{r})e^\al,z)uz^{j+1}\frac{(1+z)^{m}}{z^{2+n}}\\
	&\quad +\Res_{z}\sum_{j\geq 0} a(-j-1)Y(\al(-n_{2})\ds \al(-n_{r})e^\al,z)uz^{j+1}\frac{(1+z)^{m}}{z^{2+n}}\\
	&= \Res_{z}[a(-1), Y(\al(-n_{2})\ds \al(-n_{r})e^\al, z)]u\frac{(1+z)^{m}}{z^{2+n}}\\
	&\quad + \Res_{z}Y(\al(-n_{2})\ds \al(-n_{r})e^\al,z)a(-1)u\frac{(1+z)^{m}}{z^{2+n}}\\
	&\quad +\Res_{z}\sum_{j\geq 0} (a(-j-2)+a(-j-1))Y(\al(-n_{2})\ds \al(-n_{r})e^\al,z)uz^{j+1}\frac{(1+z)^{m}}{z^{2+n}}\\
	&\equiv \sum_{i\geq 0} \binom{-1}{i}z^{-1-i} \Res_{z} Y(a(i)\al(-n_{2})\ds \al(-n_{r})e^\al,z)u \frac{(1+z)^{m}}{z^{2+n}}\pmod{O'}\\
	&=\Res_{z} 2Y(\al(-n_{2})\ds \al(-n_{r})e^\al)u \frac{(1+z)^m}{z^{3+n}}\\
	&\quad + \sum_{i\geq 0}\sum_{s=2}^{r} \binom{-1}{i} \Res_{z} Y(\al(-n_{2})\ds[\al(i),\al(-n_{s})]\ds \al(-n_{r})e^\al,z)u \frac{(1+z)^{m}}{z^{2+n+i+1}}\\
	&\equiv \sum_{s=2}^{r}(-1)^{n_s} n_{s} \Res_{z} Y(\al(-n_{2})\ds \widehat{\al(-n_{s})}\ds \al(-n_{r})e^\al,z)u \frac{(1+z)^{m}}{z^{2+n+n_s+1}}\pmod{O'}.
\end{align*}
Denote $\al(-n_{2})\ds \widehat{\al(-n_{s})}\ds \al(-n_{r})e^\al$ by $a_{s}$. We have $m=\wt a_{s}+n_{s}$. Then by the induction hypothesis \eqref{3.17}, with $r$ replaced by $r-1$, we have 
\begin{align*}
	&\sum_{s=2}^{r}(-1)^{n_s} n_{s} \Res_{z} Y(\al(-n_{2})\ds \widehat{\al(-n_{s})}\ds \al(-n_{r})e^\al,z)u \frac{(1+z)^{m}}{z^{2+n+n_s+1}}\\
	&=\sum_{s=2}^{r} (-1)^{n_{s}} n_{s} \Res_{z} Y(a_{s},z)u(1+z)^{n_{s}}\frac{(1+z)^{\wt a_{s}}}{z^{2+n+n_s+1}}\\
	&=\sum_{s=2}^{r} \sum_{j\geq 0}\binom{n_{s}}{j} (-1)^{n_{s}} n_{s} \Res_{z} Y(a_{s},z)u\frac{(1+z)^{\wt a_{s}}}{z^{2+n_s+j+1}}\\
	&\equiv 0\pmod{O'},
\end{align*}
since $n_s+j+1\geq 1$. This proves \eqref{3.18}. Now assume that 
\begin{equation}\label{3.19}
	\Res_{z} Y(\al(-k)\al(-n_{2})\ds \al(-n_{r})e^\al,z)u \frac{(1+z)^{N+k+n_{2}+\ds+n_{r}}}{z^{2+n}}\in O',
\end{equation}
for some fixed $k\geq 1$, $n_{2},\ds,n_{r}\geq 1$, and $n\geq 0$. We want to show that 
\begin{equation}\label{3.20}
	\Res_{z} Y(\al(-k-1)\al(-n_{2})\ds \al(-n_{r})e^\al,z)u \frac{(1+z)^{N+k+1+n_{2}+\ds+n_{r}}}{z^{2+n}}\in O'.
\end{equation}
Indeed, by adopting a similar calculation as the proof of Lemma \ref{lm3.6}, we have
\begin{align*}
	&\Res_{z} Y(\al(-k-1)\al(-n_{2})\ds \al(-n_{r})e^\al,z)u \frac{(1+z)^{N+k+1+n_{2}+\ds+n_{r}}}{z^{2+n}}\\
	&=\Res_{z}\frac{1}{k}Y(L(-1)\al(-k)\al(-n_{2})\ds \al(-n_{r})e^\al,z)u \frac{(1+z)^{N+k+1+n_{2}+\ds +n_{r}}}{z^{2+n}}\\
	&\quad +\Res_{z}\frac{1}{k}Y(\al(-k)[L(-1),\al(-n_{2})\ds \al(-n_{r})]e^\al,z)u \frac{(1+z)^{N+k+1+n_{2}+...+n_{r}}}{z^{2+n}}\\
	&\quad +\Res_{z} \frac{1}{k} Y(\al(-1)\al(-k)\al(-n_{2})\ds \al(-n_r)e^\al)u\frac{(1+z)^{N+k+1+n_{2}+\ds+n_{r}}}{z^{2+n}}\\
	&=-\Res_{z}\frac{1}{k}Y(\al(-k)\al(-n_{2})\ds \al(-n_{r})e^\al,z)u \frac{d}{dz}\left( \frac{(1+z)^{N+k+1+n_{2}+\ds+n_{r}}}{z^{2+n}}\right)\\
	&\quad +\sum_{s=2}^{r}\Res_{z}\frac{n_{s}}{k}Y(\al(-k)\al(-n_{2})\ds \al(-n_{s}-1)\ds \al(-n_{r})e^\al,z)u \frac{(1+z)^{N+k+\ds(1+n_{s})\ds+n_{r}}}{z^{2+n}}\numberthis\label{3.21}\\
	&\quad +\Res_{z} \frac{1}{k} Y(\al(-1)\al(-k)\al(-n_{2})\ds \al(-n_r)e^\al)u\frac{(1+z)^{N+k+1+n_{2}+\ds+n_{r}}}{z^{2+n}}\\
	&=-\Res_{z}\frac{1}{k}(N+k+1+\ds+n_{r})Y(\al(-k)\al(-n_{2})\ds \al(-n_{r})e^\al,z)u \frac{(1+z)^{N+k+n_{2}+\ds+n_{r}}}{z^{2+n}}\\
	&\quad +\Res_{z}\frac{2+n}{k}Y(\al(-k)\al(-n_{2})\ds \al(-n_{r})e^\al,z)u \frac{(1+z)^{N+k+n_{2}+\ds+n_{r}}(1+z)}{z^{2+n+1}}\\
	&\quad +\sum_{s=2}^{r}\Res_{z}\frac{n_{s}}{k}Y(\al(-k)\al(-n_{2})\ds \al(-n_{s}-1)\ds \al(-n_{r})e^\al,z)u \frac{(1+z)^{N+k+\ds(1+n_{s})\ds+n_{r}}}{z^{2+n}}\\
	&\quad +\Res_{z} \frac{1}{k} Y(\al(-1)\al(-k)\al(-n_{2})\ds \al(-n_r)e^\al)u\frac{(1+z)^{N+k+1+n_{2}+\ds+n_{r}}}{z^{2+n}}\\
	&\equiv 0+\Res_{z} \frac{1}{k} Y(\al(-1)\al(-k)\al(-n_{2})\ds \al(-n_r)e^\al)u\frac{(1+z)^{N+k+1+n_{2}+\ds+n_{r}}}{z^{2+n}}\pmod{O'},
\end{align*}
where the congruences follow from the induction (on $k\geq 1$) hypothesis \eqref{3.19}. Moreover, by adopting a similar argument as our previous proof of \eqref{3.18} under the given assumption \eqref{3.17}, we can show the following fact under the given assumption \eqref{3.19}: 
$$\Res_{z} \frac{1}{k} Y(\al(-1)\al(-k)\al(-n_{2})\ds \al(-n_r)e^\al)u\frac{(1+z)^{N+k+1+n_{2}+\ds+n_{r}}}{z^{2+n}}\in O',$$
Thus \eqref{3.20} is true, in view of our calculation \eqref{3.21}. Now the induction step for $k\geq 1$ and the induction step for the length $r\geq 1$ of $a\in M_{\hat{\h}}(1,\al)$ are both complete. 
\end{proof}

Now we have finished the proof of \eqref{eq:main} and Theorem~\ref{thm3.8}.

\subsection{The structure of $A(V_B)$}
In this subsection, we use Theorem~\ref{thm3.8} to prove the claimed isomorphism \eqref{eq:6.1} for the Zhu's algebra $A(V_B)$. Recall that $V_B=V_{\Z_{\geq 0}\al}=\bigoplus_{k=0}^\infty 	M_{\hat{\h}}(1,k\al)$ is the Borel-type subVOA of the lattice VOA $V_{\Z\al}$, with $(\al|\al)=2N$.


\begin{prop}\label{prop3.1}
There exists an epimorphism of associative algebras: 
\begin{equation}\label{3.6}
	F:\C\<x,y\>/\<y^2,yx+Ny,xy-Ny\>\ra A(V_B),
\end{equation}
such that $F(x)=[\al(-1)\vac]$ and $F(y)=[e^\al]$. 
\end{prop}
\begin{proof}
By the definition of $Y(e^\al,z)$ in \eqref{2.4}, for any $n\geq 0$, we have
\begin{equation}\label{3.7}
	e^\al_{n}e^\al=0, \quad e^\al_{-1}e^\al=\ds=e^\al_{-2N}e^\al=0,\quad n\geq 0,\quad \mathrm{and}\quad e^\al_{-2N-1}e^\al=e^{2\al}.
\end{equation}
Since $\wt e^\al=N$, by \eqref{3.2} and \eqref{3.7}, we have
$e^\al \ast e^\al=\sum_{j\geq 0} \binom{N}{j} e^\al_{j-1}e^\al=0$. Hence $[e^\al]\ast [e^\al]=0$ in $A(V_B)$. By \eqref{3.5}, we have $$\al(-n-2)u+\al(-n-1)u \equiv 0\pmod{O(V_B)},$$
and $[\al(-1)u]=[u]\ast [\al(-1)\vac]$, for all $n\geq 0$ and $u\in V$. Thus
\begin{align*}
	[\al(-n_{1}-1)\al(-n_{2}-1)\dots \al(-n_{k}-1)u]=(-1)^{n_1+\dots +n_k} [u]\ast [\al(-1)\vac]\ast \dots \ast [\al(-1)\vac],
\end{align*}
for any $n_{1},\dots ,n_k\geq 0$ and $u\in V$. This shows $A(V_B)$ is generated by $[\al(-1)\vac]$ and $[e^{m\al}]$, for all $m\geq  1$. We claim that $[e^{k\al}]=0$ for any $k\geq 2$. 

Indeed, for $m\geq 1$, since we have $e^\al_{-2Nm-1}e^{m\al}=e^{(m+1)\al}$, $e^{\al}_{-n}e^{m\al}=0$ for any $n\leq 2Nm$, and $2Nm+1\geq 2$, then it follows from \eqref{3.5} that for any $m\geq 1$, 
\begin{align*}e^{(m+1)\al}&=e^{\al}_{-2Nm-1}e^{m\al}+\binom{N}{1}e^\al_{-2Nm}e^{m\al}+\ds +\binom{N}{N}e^\al_{-2Nm-1+N}e^{m\al}\\
	&=\Res_{z} Y(e^\al,z)e^{m\al}\frac{(1+z)^{N}}{z^{2Nm+1}} \equiv 0\pmod{O(V_B)}. 
\end{align*}
Hence $[e^{k\al}]=0$ in $A(V_B)$ for all $k\geq 2$, and $A(V_B)$ is generated by $[\al(-1)\vac]$ and $[e^\al]$. Then we have an epimorphism $F:\C\<x,y\>\ra A(V_B),$ such that $F(x)=[\al(-1)\vac]$ and $F(y)=[e^\al]$. Moreover, by \eqref{3.1} and the definition of $Y(e^\al,z)$ in \eqref{2.4}, we have 
\begin{align*}
	e^\al\circ \vac&=e^\al_{-2}\vac+\binom{N}{1}e^\al_{-1}\vac+\sum_{j\geq 2}\binom{N}{j}e^\al_{j-2} \vac =\Res_{z} z^{-2}\exp(-\sum_{n<0}\frac{\al(n)}{n}z^{-n})e^\al+Ne^\al+0\\
	&= \al(-1)e^\al+Ne^\al\equiv 0\pmod{ O(V_B)},
\end{align*}
and so $[e^\al]\ast [\al(-1)\vac]+N[e^\al]=0$ in $A(V_B)$. By \eqref{3.4}, we also have
$$[\al(-1)\vac]\ast [e^\al]-[e^\al]\ast [\al(-1)\vac]=[\Res_{z}Y(\al(-1)\vac,z)e^\al ]=[\al(0)e^\al]=2N[e^\al],$$
and so $[\al(-1)\vac]\ast [e^\al]-N[e^\al]=0$ in $A(V_B)$. Therefore, the epimorphism $F:\C\<x,y\>\ra A(V_B)$ factors through $\C\<x,y\>/\<y^2,yx+Ny,xy-Ny\>$, and induces an epimorphism $F$ in \eqref{3.6}. 
\end{proof}

In fact, the proof of Proposition~\ref{prop3.1} does not use the concrete spanning set of $O(V_B)$. But to show that $F$ defined by \eqref{3.6} is an isomorphism, we have to use Theorem~\ref{thm3.8}.

\begin{coro}\label{coro3.9}
The epimorphism $F$ given by \eqref{3.6} is an isomorphism of associative algebras. In particular, we have $A(V_B)\cong \C[x]\op \C y$, with 
\begin{equation}\label{3.23}
	y^2=0,\quad yx=-Ny,\quad xy=Ny.
\end{equation}
\end{coro}
\begin{proof}
We construct an inverse map of $F$ in \eqref{3.6}. Let
\begin{align*}
	G: V_B\ra\C\<x,y\>/&\<y^2,yx+Ny,xy-Ny\>, \\
	\al(-n_{1}-1)\dots \al(-n_{r}-1)\vac &\mapsto (-1)^{n_{1}+\dots n_r} x^r,\numberthis \label{3.24}\\
	\al(-n_1-1)\dots \al(-n_r-1)e^\al& \mapsto (-1)^{n_1+\dots n_r} yx^r=(-1)^{r+n_1+\ds n_r} y,\\
	M_{\hat{\h}}(1,k\al)&\mapsto 0, \quad k\geq 2,
\end{align*}
where $r\geq 0$ and $n_1,\ds ,n_r\geq 0$, and we use the same symbols $x$ and $y$ to denote their image in the quotient space. Note that $G$ is well-defined since $V_B=\bigoplus_{k\geq 0} M_{\hat{\h}}(1,k\al)$, and $	\al(-n_{1}-1)\dots \al(-n_{r}-1)\vac$ and $	\al(-n_1-1)\dots \al(-n_r-1)e^\al$ are basis elements of $M_{\hat{\h}}(1,0)$ and $M_{\hat{\h}}(1,\al)$, respectively. We claim that $G(O(V_B))=0$. 

Indeed, it suffices to show that $G$ vanishes on the spanning elements of $O(V_B)$ in \eqref{3.22}. By the definition \eqref{3.24} of $G$, we already have $G(	M_{\hat{\h}}(1,k\al))=0$ for any $k\geq 2$. In particular, we have $G(\al(-n-2)u+\al(-n-1)u)=G(\al(-1)v+Nv)=0$ if $u,v\in 	M_{\hat{\h}}(1,k\al)$ for some $k\geq 2$.

If $u=	\al(-n_{1}-1)\dots \al(-n_{r}-1)\vac\in M_{\hat{\h}}(1,0)$, then it follows from \eqref{3.24} that
\begin{align*}
	&G(\al(-n-2)u+\al(-n-1)u)\\
	&=G(\al(-n-2)\al(-n_{1}-1)\dots \al(-n_{r}-1)\vac)+G(\al(-n-1)\al(-n_{1}-1)\dots \vac)\\
	&=(-1)^{n+1+n_1+\ds+n_r} x^{r+1}+(-1)^{n+n_1+\ds +n_r} x^{r+1}=0.
\end{align*}

If $u=	\al(-n_1-1)\dots \al(-n_r-1)e^\al\in M_{\hat{\h}}(1,\al)$, then it follows from \eqref{3.24} that
\begin{align*}
	&G(\al(-n-1)u+\al(-n-1)u)\\
	&=G(\al(-n-2)	\al(-n_1-1)\dots \al(-n_r-1)e^\al)+G(\al(-n-1)	\al(-n_1-1)\dots e^\al)\\
	&=(-1)^{n+1+n_1+\ds+n_r}y x^{r+1}+(-1)^{n+n_1+\ds +n_r} yx^{r+1}=0.
\end{align*}
Thus, $G(\al(-n-2)u+\al(-n-1)u)=0$ for any $u\in V_B$. Finally, if $v=\al(-n_{1}-1)\dots \al(-n_{r}-1)e^\al\in M_{\hat{\h}}(1,\al)$, it again follows from \eqref{3.24} that
\begin{align*}
	&G(\al(-1)v+Nv)\\
	&=G(\al(-1)\al(-n_{1}-1)\dots \al(-n_{r}-1)e^\al)+NG(\al(-n_{1}-1)\dots \al(-n_{r}-1)e^\al)\\
	&=(-1)^{n_1+\ds n_r} yx^{r+1}+(-1)^{n_1+\ds n_r} Nyx^r\\
	&=(-1)^{n_1+\ds n_r} (yx+Ny)x^r=0,
\end{align*}
as $yx+Ny=0$. Thus, $G$ given by \eqref{3.24} induces a linear map 
\begin{align*}
	G:A(V_B)=V_B/O(V_B)&\ra \C\<x,y\>/\<y^2,yx+Ny,xy-Ny\>,\quad \mathrm{such\ that}\\
	G([\al(-n_{1}-1)\dots \al(-n_{r}-1)\vac])&=G((-1)^{n_1+\ds +n_r}[\al(-1)\vac]^{r})=(-1)^{n_1+\ds +n_r}x^r,\numberthis \label{3.25}\\
	G([\al(-n_{1}-1)\dots \al(-n_{r}-1)e^\al])&=G((-1)^{n_1+\ds +n_r}[e^\al]\ast[\al(-1)\vac]^{r})=(-1)^{n_1+\ds +n_r}yx^r,
\end{align*}
for any $r\geq 0$, $n_1,\ds ,n_r\geq 0$, and $k\geq 2$. Since $A(V_B)$ is spanned by elements of the form $[\al(-n_{1}-1)\dots \al(-n_{r}-1)\vac]$ and $[\al(-n_{1}-1)\dots \al(-n_{r}-1)e^\al]$ because of \eqref{3.22}, it is clear that $G\circ F=\Id$ and $F\circ G=\Id$, in view of \eqref{3.6} and \eqref{3.25}. 
\end{proof}

\subsection{Applications of the main theorem}\label{Sec:7.3}
As an immediate application of the structural theorem for $A(V_B)$, we can classify the irreducible modules over $V_B$ and compute the fusion rules among them. We can also present a spanning set for $O(V_{A_1})$, which was unknown.

\subsubsection{Classification of irreducible modules over $V_B$}

\begin{lm}\label{lm5.2.10}
If $U\neq 0$ is an irreducible module over $A(V_{B})\cong \C[x]\op \C y$, then we must have $y.U=0$, and $U\cong \C e^\la$ for some $\la\in \h=\C\al$, with $x.e^\la=(\al|\la)e^\la$. 
\end{lm}
\begin{proof}
By \eqref{3.23}, $\C y$ is an ideal of $A(V_{\Z_{\geq 0}\al})$. Then $y.U\leq U$ is an $A(V_B)$-submodule, and so $y.U$ is either $U$ or $0$. If $y.U=U$, then we have $0=y^2.U=y.U=U$, a contradiction. Thus, $y.U=0$ and $U$ is an irreducible module over $\C[x]$. We have $U\cong \C[x]/\mathfrak{m}$, for some maximal ideal $\mathfrak{m}$ of $\C[x]$. By Hilbert's Nullstellensatz, we have $\mathfrak{m}=\<x-\mu\>$ for some $\mu \in \C$. Choose $\la\in \h$ so that $(\al|\la)=\mu$. Then $U\cong \C[x]/\<x-(\al|\la)\>\cong \C e^\la$, with $x.e^\la=(\al|\la) e^\la$. 
\end{proof}

\begin{lm}\label{lm3.13}
For any irreducible module $W=M_{\hat{\h}}(1,\la)$ over the Heisenberg VOA $M_{\hat{\h}}(1,0)$, $W$ is also an irreducible module over the Borel-type subVOA $V_{B}$, where $Y_W:V_{B}\ra \End(W)[[z,z^{-1}]]$ satisfies $Y_W(a,z)=0$, for any $a\in M_{\hat{\h}}(1,n\al)$ and $n\geq 1$, and $Y_W|_{M_{\hat{\h}}(1,0)}$ is given by the action of the Heisenberg VOA $M_{\hat{\h}}(1,0)$. 
\end{lm}
\begin{proof}
By \eqref{2.8}, for any $a\in M_{\hat{\h}}(1,n\al)$ and $b\in M_{\hat{\h}}(1,m\al)$, where $m,n\geq 0$, $Y_W(a,z)b$ is either $0$ or given by the Heisenberg module vertex operator. Hence $(W,Y_W)$ is a well-defined module over the Borel-type subVOA $V_{B}$. It is clear that $W$ is irreducible. 
\end{proof}

\begin{thm}\label{thm3.14}
$\Sigma=\left\{(W=M_{\hat{\h}}(1,\la),Y_W): \la\in \h=\C\al\right\}$, with $Y_W$ defined by Lemma~\ref{lm3.13}, is a complete list of irreducible modules over the rank-one Borel-type subVOA $V_{B}$. 

Moreover, the fusion rule between the irreducible $V_{B}$-modules $M_{\hat{\h}}(1,\la), M_{\hat{\h}}(1,\mu)$, and $M_{\hat{\h}}(1,\gamma)$ is the same as the fusion rule  of these modules as modules over the Heisenberg VOA. In other words, $N\fusion{M_{\hat{\h}}(1,\la)}{M_{\hat{\h}}(1,\mu)}{M_{\hat{\h}}(1,\gamma)}\cong \delta_{\la+\mu,\gamma}$. 
\end{thm}
\begin{proof}
Given a module $(W=M_{\hat{\h}}(1,\la),Y_W)$ in $\Sigma$, the bottom degree $W(0)=\C e^\la$ is an $A(V_{B})$-module, with the actions of $x=[\al(-1)\vac]$ and $y=[e^\al]$ given by 
$$x.e^\al=o(\al(-1)\vac)e^\la=(\al|\la)e^\la, \quad  y.e^\la=o(e^\al)e^\la=\Res_{z}z^{N-1} Y_W(e^\al,z)e^\la=0,$$
see \cite{FZ,Z}. By Lemma~\ref{lm3.13}, such $A(V_{B})$-modules $W(0)$, with $W$ varies in $\Sigma$, are all the irreducible modules over $A(V_B)$ up to isomorphism. Then by Theorem 2.2.2 in \cite{Z}, $\Sigma$ is the complete list of irreducible modules over $V_{B}$. Finally, note that any intertwining operator between modules over the Heisenberg VOA $I\in I\fusion{M_{\hat{\h}}(1,\la)}{M_{\hat{\h}}(1,\mu)}{M_{\hat{\h}}(1,\gamma)}$ can be naturally lifted up to an intertwining operator $\tilde{I}$ of $V_{B}$, since the Jacobi identity of $I$ is 
\begin{equation}\label{5.2.26}
	\begin{aligned}
		&z_{0}^{-1}\delta\left(\frac{z_{1}-z_{2}}{z_{0}}\right) Y_{W^3}(a,z_{1})I(v,z_{2})u-z_{0}^{-1}\delta\left(\frac{-z_{2}+z_{1}}{z_{0}}\right)I(v,z_{2})Y_{W^2}(a,z_{1})u\\
		&=z_{2}^{-1}\delta \left(\frac{z_{1}-z_{0}}{z_{2}}\right) I(Y_{W^1}(a,z_{0})v,z_{2})u,
	\end{aligned}
\end{equation}
and $Y_{W^i}(a,z)=0$ for $i=1,2,3$, if $a\in M_{\hat{\h}}(1,n\al)$ with $n\geq 1$. Therefore, we can replace $I$ in \eqref{5.2.26} by the intertwining operator $\tilde{I}$ of $V_{B}$. Conversely, we can also restrict any intertwining operator $\Y\in I\fusion{M_{\hat{\h}}(1,\la)}{M_{\hat{\h}}(1,\mu)}{M_{\hat{\h}}(1,\gamma)}$ of $V_{B}$ to an intertwining operator $\Y$ of the same type between modules over the Heisenberg VOA. Therefore, the fusion rules of the Borel-type subVOA $V_{B}$  is the same as the fusion rules between modules over the Heisenberg VOA $M_{\hat{\h}}(1,0)$. 
\end{proof}

\begin{remark}
Theorem~\ref{thm3.14} is also parallel to the semisimple Lie algebra case. Note that a Borel subalgebra $\mathfrak{b}=\mathfrak{n}_+\op \h$ of a semisimple Lie algebra $\g$ has the same irreducible modules as its Cartan part $\h$.
\end{remark}

\subsubsection{A spanning set of $O(A_1)$}
For the $A_1$-root lattice $A_1=\Z\al$ with $(\al|\al)=2$, the following Corollary presents a spanning set of  $O(V_{A_1})$.
\begin{coro}\label{coro3.10}
Let $L=A_1=\Z\al$. Then $O(V_{A_1})$ is spanned by the following elements: 
\begin{equation}\label{3.29}
	\begin{cases}
		&\al(-n-2)u+\al(-n-1)u, \qquad u\in V_{A_1},\ \mathrm{and}\ n\geq 0,\\
		&\pm \al(-1)v+v,\qquad v\in  M_{\hat{\h}}(1, \pm\al),\\
		&M_{\hat{\h}}(1,\pm k\al),\qquad k\geq 2.\\
		& \al(-1)^3w-\al(-1)w,\qquad w\in M_{\hat{\h}}(1,0).
\end{cases}\end{equation}
\end{coro}
\begin{proof}
Let $O''$ be the subspace of $V_{A_1}$ spanned by the elements \eqref{3.29}. Since both $O(V_{\Z_{\geq 0}\al})$ and $O(V_{\Z_{\geq 0}(-\al)})$ are contained in $O(V_{A_1})$, 
by Lemma~\ref{lm3.2} and \ref{lm3.3}, we have $\al(-n-2)u+\al(-n-1)u\in O(V_{A_1})$ for any $u\in V_{A_1}$, $\pm \al(-1)v+v\in  O(V_{A_1})$ for any $v\in  M_{\hat{\h}}(1,\pm\al)$, and $M_{\hat{\h}}(1,\pm k\al)\subset O(V_{A_1})$ for all $k\geq 2$. Moreover, since $\wt e^{\al}=2$ and $(\al|-\al)=-2$, we have 
\begin{align*}
	e^{\al}\circ e^{-\al}&=e^{\al}_{-2}e^{-\al}+ e^{\al}_{-1}e^{-\al}\\
	&=\frac{1}{3}\al(-3)\vac+\frac{1}{4}(\al(-2)\al(-1)+\al(-1)\al(-2))\vac+ \frac{1}{3!} \al(-1)^3\vac+\frac{1}{2}\al(-2)\vac+\frac{1}{2!}\al(-1)^2\vac\\
	&\equiv \frac{1}{3}\al(-1)\vac-\frac{1}{2}\al(-1)^2\vac+\frac{1}{6}\al(-1)^3\vac-\frac{1}{2}\al(-1)\vac+\frac{1}{2}\al(-1)^2\vac\\
	&\equiv \frac{1}{6}(\al(-1)^3\vac-\al(-1)\vac)\pmod{O(V_{A_1})}.
\end{align*}
Hence $\al(-1)^3w-\al(-1)w\in O(V_{A_1})$ for any $w\in M_{\hat{\h}}(1,0)$, in view of the proof of Lemma~\ref{lm3.3}. This shows $O''\ssq O(V_{A_1})$. Conversely, with a similar argument as Proposition~\ref{prop3.1} and Corollary~\ref{coro3.9}, we can easily show that 
$$V_{A_1}/O''\cong \C z\op \C[x]/\<x^3-x\>\op \C y,\quad [e^{-\al}]\mapsto z,\ [\al(-1)\vac]\mapsto x,\ [e^{\al}]\mapsto y.$$
In particular, we have $\dim V_{A_1}/O''=5$. On the other hand, we have $\dim V/O(V_{A_1})=5$, since $V_{A_1}$ has two irreducible modules $V_{\Z\al}$ and $V_{\Z\al+\frac{1}{2}\al}$, whose bottom degrees are of dimensions $1$ and $2$, respectively, see \cite[Theorem 3.1]{D}. Using the one-to-one correspondence between irreducible $A(V_{A_1})$-modules and irreducible $V_{A_1}$-modules (cf. \cite[Theorem 2.2.2]{Z}), we have $O''=O(V_{A_1})$, as there is an epimorphism $V_{A_1}/O''\ra A(V_{A_1})$. 
\end{proof}


\section{Acknowledgements}
I'm grateful to Professors Angela Gibney and Daniel Krashen for helpful discussions and many useful comments on this paper. This paper is an expansion of a Chapter in my Ph.D. dissertation submitted to the University of California, Santa Cruz. I would like to thank my Ph.D. advisor, Professor Chongying Dong, for his guidance and advice on my research in vertex operator algebras.




\end{document}